\documentclass[12pt,reqno,notitlepage]{amsart}
\usepackage{amsmath,amsfonts,amsthm,amssymb}
\usepackage{enumerate}
\usepackage{indentfirst}
\usepackage[dvips]{graphicx}
\usepackage[all]{xy}
\usepackage{stackrel}
\usepackage{marginnote}

\usepackage[centering, includeheadfoot, hmargin=1.0in, tmargin=1.0in, 
bmargin=1in, headheight=6pt]{geometry}
\usepackage[latin1]{inputenc}
\usepackage[T1]{fontenc}
\usepackage{verbatim}
\usepackage{color}
\usepackage[normalem]{ulem}

\usepackage{hyperref}  
\hypersetup{
pdfborder={0 0 0}, 
colorlinks=true, 
citecolor=blue,
linktoc=page,
pdfauthor={Thiago Fassarella and Frank Loray}, 
pdftitle={Higgs bundles}
}
\renewcommand{\ref}{\hyperref}

\newtheorem{thm}{Theorem}[section]

\newtheorem{prop}[thm]{Proposition}

\newtheorem{lemma}[thm]{Lemma}

\newtheorem{cor}[thm]{Corollary}
\newtheorem{remark}[thm]{Remark}

\theoremstyle{definition}
\newtheorem{definition}[thm]{Definition}

\renewcommand{\P}{\mathbb{P}}
\newcommand{\C}{\mathbb{C}}

\numberwithin{equation}{section}

\begin{document}
\title{Moduli of Higgs bundles over the five punctured sphere}
\author[T. Fassarella]{Thiago Fassarella}
\address{\color{black}Universidade Federal Fluminense, Instituto de Matem\'atica e Estat\'istica.
Rua Alexandre Moura 8, S\~ao Domingos, 24210-200 Niter\'oi RJ,
Brazil.}
\email{\color{black}tfassarella@id.uff.br}

\author[F. Loray]{Frank Loray}
\address{Univ Rennes, CNRS, IRMAR, UMR 6625, F-35000 Rennes, France.}
\email{frank.loray@univ-rennes1.fr}

\subjclass[2010]{Primary 34M55; Secondary 14D20, 32G20, 32G34.}
\keywords{Higgs bundles, parabolic structure, elliptic curve, Hitchin fibration, spectral curve. }
\thanks{The second author is supported by CNRS and Centre Henri Lebesgue, program ANR-11-LABX-0020-0. 
 The authors also thank  Brazilian-French Network in Mathematics and CAPES-COFECUB  932/19.}
 \date{\today}

\begin{abstract}
We look at rank two parabolic Higgs bundles over the projective line minus five points which are semistable with respect to a weight vector $\mu\in[0,1]^5$. The moduli space corresponding to the central weight  $\mu_c=(\frac{1}{2}, \dots, \frac{1}{2})$ is studied in details and all singular fibers of the Hitchin map are described, including the nilpotent cone.  After giving a description of fixed points of the $\C^*$-action we obtain a proof of Simpson's foliation conjecture in this case. For each $n\ge 5$, we remark that there is a weight vector so that the foliation conjecture in the moduli space of rank two logarithmic connections over the projective line minus $n$ points is false. 
\end{abstract}

\maketitle


\section{Introduction}

Let $\P^1$ be the complex projective line and let $\Lambda=0+1+\lambda+t+\infty$ be a  divisor on it,  supported on five distinct points. We refer to $\Lambda$ as divisor of parabolic points. Let $\omega_{\P^1}(\Lambda)$ be the sheaf of 1-forms with simple poles on $\Lambda$. We study moduli spaces   of traceless semistable  parabolic Higgs bundles over $(\P^1, \Lambda)$. These moduli spaces parametrize triples $(E, {l}, \theta)$,  where $E$ is a rank two vector bundle on  $\P^1$, of degree zero,  with the additional data  ${l} = \{l_i\}$, a one dimensional subspace $l_i$,  on the fiber over each parabolic point, and $\theta: E\to E\otimes \omega_{\P^1}(\Lambda)$ consists of a traceless  homomorphism, which is called a Higgs field.  The construction of these moduli spaces depends on the choice of a weight vector $\mu\in [0,1]^5$, which determines the stability condition.

The nonabelian Hodge correspondence \cite{H87, Sim90, Sim92}, gives an identification between Higgs bundles and local systems, then our moduli spaces are in correspondence with the character variety, which parametrizes local systems on $\P^1$ with given conjugacy classes at the punctures.  This last is a topological object, while in  the former case it depends on the algebraic structure of the punctured line. This leads to the investigation of algebraic invariants, such as the Hitchin fibration, which is an important object with recent deep results, see for example \cite{HT, Ngo1, CHM12}.  The Hitchin map is defined on this moduli spaces, it sends a Higgs field to its determinant, which is a quadratic differential. On the one hand, this map is known to be an algebraically integrable system, i.e. its general fiber is Lagrangian and isomorphic to an  Abelian variety. On the  other hand,  singular fibers are difficult to deal with. In \cite{Hor, Hor2, GO13}, singular fibers in the moduli space of twisted pairs were studied.   In this paper we determine all singular fibers of the Hitchin map in the specific parabolic case $(\P^1, \Lambda)$.   

We also describe the locus of fixed points with respect to the $\C^*$-action given by multiplication on the Higgs field. This motivates to investigate the foliation conjecture \cite[Question 7.4]{S08} on the moduli space of rank two logarithmic connections with generic residues. The parabolic version of the foliation conjecture has been proved in \cite[Corollaries 5.7 and 6.2]{LSS} for the moduli space of connections over $\P^1$ minus four points (when the weight vector is generic),  and recently \cite{HHZ} deals with the case $\P^1$ minus five points by assuming  the weight vector $\mu$ satisfies $\sum\mu_i<1$. Since this last case lies in the unstable zone (any  parabolic connection has $\mu$-unstable parabolic vector bundle) we turn our attention to  the stable zone. After determining the locus of fixed points of the $\C^*$-action, we obtain a proof of the foliation conjecture in the case $\P^1$ minus five points with the central weight vector $\mu_c=\left(\frac{1}{2}, \dots, \frac{1}{2} \right)$. We also remark that for all $n\ge 5$ there is a weight vector (in the stable zone),  such that the foliation conjecture in the case $\P^1$ minus $n$ points is false. 
  

In our context, every Higgs field having nonvanishing determinant is irreducible, i.e. it does not have any invariant line subbundle, then it is stable for any choice of weight vector. The moduli space $\mathcal H$ associated to the central weight $\mu_c=\left(\frac{1}{2}, \dots, \frac{1}{2} \right)$ is particularly interesting,  indeed it is a smooth quasi-projective variety of dimension four and its automorphism group admits a modular realization of $\big(\mathbb Z/2\mathbb Z\big)^4$ as a subgroup. This subgroup, denoted here by ${\bf El}$, consists  of elementary transformations $elem_I$,  for each subset  $I\subset \{0,1,\lambda, t, \infty\}$ of  even cardinality. 

We shall consider only Higgs fields which are nilpotent with respect to the parabolic direction.  This implies that our moduli space  $\mathcal H$ contains an open dense subset $\mathcal U$ isomorphic to the cotangent bundle $T^*\mathcal S$, where $\mathcal S$ denotes the moduli space of parabolic vector bundles.  It is well known that $\mathcal S$ is a del Pezzo surface of degree four, see \cite{B, LS}, its automorphism group has order $16$ and coincides with the group ${\bf El}$ of elementary transformations \cite{LS, AFKM21}. There are  16 rational curves $\zeta_i$ with $(-1)$-self intersection on this surface, we denote by $\Sigma$ the union of them.

The main goal of this paper is to determining all singular fibers of the Hitchin map. The most complicated one is the nilpotent cone $\mathcal N$, consisting of Higgs fields having vanishing determinant. In order to describe it,  let us  consider the forgetful map 
$\frak{for} : \mathcal H \dashrightarrow \mathcal S$,
which forgets the Higgs field.  Note that, since $\mathcal S$ admits an embedding in $\mathcal H$, by taking the Higgs field to be zero, it gives one component of  $\mathcal N$. Our first goal is the following result.

\begin{thm}\label{thmintro:nilpP1}
The nilpotent cone of $\mathcal H$ has exactly $17$ irreducible components
\[
\mathcal N = \mathcal S\cup_{i=1}^{16}\mathcal N_i
\]
where $\frak{for}(\mathcal N_i) = \zeta_i$.  See Figure~\ref{nilpcone}.
\end{thm} 
 
This is Theorem~\ref{thm:nilpP1} in the main text.  
 
Before describing the remaining singular fibers, let us briefly introduce the spectral curve.   The Hitchin basis, formed by quadratic differentials, is two dimensional. The locus of singular spectral curves is a union of five lines. For instance, the general spectral curve $X_s$ is a smooth curve of genus two, branched over six points  $0,1,\lambda, t, \infty, \rho$ of $\P^1$ and the corresponding Hitchin fiber is isomorphic to the Picard variety ${\rm Pic}^{3}(X_s)$ which parametrizes degree 3 line bundles on $X_s$. 

A singular spectral curve occurs when the sixth point $\rho$ coincides with one of the five other points.  This leads to a nodal curve $X_s$ of genus $2$, whose desingularization $\tilde{X_s}$ is an elliptic curve branched over 
\[
\{0,1,\lambda, t, \infty\}\setminus\{\rho\} 
\]
and $X_s$ can be obtained identifying two points $w_{\rho}^+$ and $w_{\rho}^-$ of $\tilde{X_s}$.  Let us mention that the compactified Jacobian $\overline{{\rm Pic}}^0(X_s)$, which parametrizes isomorphism classes of torsion free sheaves of rank one and degree zero on $X_s$,  is obtained identifying the $0$-section with the $\infty$-section of the $\P^1$-bundle 
\begin{eqnarray}\label{compJac}
{\bf F} = \P(\mathcal O_{\tilde{X_s}}(w_{\rho}^+)\oplus\mathcal O_{\tilde{X_s}}(w_{\rho}^-))
\end{eqnarray}
via the translation $\mathcal O_{\tilde{X_s}}(w_{\rho}^+ - w_{\rho}^-)$ (cf. \cite[p. 83]{OS79}). See Figure~\ref{p1bundle}.

We now describe the remaining singular Hitchin fibers. For this, we consider the map 
\[
f: \mathcal H\setminus \mathcal N\to \mathcal H^{pairs}
\] 
to the moduli space of pairs $(E, \theta)$, which forgets the parabolic direction. This map consists of a blowing-up of the locus formed by Higgs fields which are holomorphic at some point $\rho\in\{0, 1, \lambda, \infty\}$ (Lemma~\ref{lemma:forgetful}). Now, let  ${\det}^{-1}(s)$ denote a singular Hitchin fiber, $s\neq 0$, coming from a singular spectral curve $X_s$ which has a node at $\rho$. We find that ${\det}^{-1}(s)$ has two irreducible components ${\bf F}_{hol}$ and ${\bf F}_{app}$, which are isomorphic to ${\bf F}$. The first one parametrizes Higgs fields which are holomorphic at $\rho$, it is contracted by $f$, the second is formed by Higgs fields which are apparent with respect to the parabolic direction over $\rho$. In addition, the restriction of $f$ to ${\bf F}_{app}$ gives a desingularization of the compactified Jacobian $\overline{\rm Pic}^3(X_s)$. This leads to the following result, which corresponds to Theorem~\ref{thm:sfiber}. 

\begin{thm}\label{thmintro:sfiber}
Assume that the spectral curve $X_s$ has a  nodal singularity at $\rho\in \{0,1,\lambda, t, \infty\}$. The corresponding singular fiber ${\det}^{-1}(s)$ of the Hitchin map  has two irreducible components
\[
{\det}^{-1} (s) = {\bf F}_{hol}\cup{\bf F}_{app}
\]
which are isomorphic via any elementary transformation
\[
(elem_I)|_{{\bf F}_{hol}} : {\bf F}_{hol} \to {\bf F}_{app}
\]
where $I \subset \{0,1,\lambda, t, \infty\}$  contains $\rho$ and has even cardinality. Moreover:
\begin{enumerate}
\item Each component is a desingularization of $\overline{\rm Pic}^3(X_s)$, then isomorphic to ${\bf F}$, and the structure of $\P^1$-bundle in ${\bf F}_{hol}$ is given by 
\[
f|_{{\bf F}_{hol}}: {\bf F}_{hol}\to \tilde{X_s}\simeq \overline{\rm Pic}^3(X_s)\setminus  {\rm Pic}^3(X_s).
\] 
\item The map $f|_{{\bf F}_{app}}: {\bf F}_{app}\to  \overline{\rm Pic}^3(X_s)$ is a desingularization map. See Figure~\ref{Fapp}.
\item The intersection ${\bf F}_{hol}\cap{\bf F}_{app}$ is the union of the $0$-section and the $\infty$-section of ${\bf F}_{hol}$. See Figure~\ref{interComp}.
\end{enumerate}
\end{thm}

In particular, we find that each component of the singular fiber $\det^{-1}(s)$ is a decomposable $\P^1$-bundle over an elliptic curve, it consists of the desingularization of  the compactified Jacobian of the corresponding nodal spectral curve. The  whole fiber topologically looks like an elliptic curve times a degenerate elliptic curve (two copies of $\P^1$ meeting in two points), but a suitable translation must be considered, see Remark~\ref{rkm:multi2}.   This confirms the guess of C.T. Simpson \cite[Discussion]{S18}.  

In the last section of the paper we deal with the moduli space  $\mathcal C^{\nu}(\P^1, \Lambda_n)$, $n\ge 5$,  of $SL_2$ logarithmic connections over $\P^1$. Here,  $\Lambda_n=t_1+\cdots+t_n$ intends to be the polar divisor, which is supported on $n$ distinct points, and $\nu = (\nu_1,\dots, \nu_n)\in\C^n$ is a  prescribed eigenvalue vector. For each weight vector $\mu$ and for $(E,\nabla,{l})\in\mathcal C^{\nu}(\P^1, \Lambda_n)$ there exists a unique limit $ \lim_{c\to 0}  c\cdot(E, \nabla, {l})$   in the moduli space of $\mu$-semistable parabolic Higgs bundles.  This gives an equivalence relation on $\mathcal C^{\nu}(\P^1, \Lambda_n)$ by assuming that two points  are equivalent if their limits are the same.  The foliation conjecture \cite[Question 7.4]{S08}, in this case, predicts that this decomposition  is  a Lagrangian (regular) foliation $\mathcal F_{\mu}$. We obtain the following result, which corresponds to Proposition~\ref{prop:foliation} and Theorem~\ref{thm:foliation}. 

\begin{thm} For the moduli space $\mathcal C^{\nu}(\P^1, \Lambda_n)$ we have:
\begin{enumerate}
\item For each $n\ge 5$ there is a weight vector $\mu$ such that the foliation conjecture is false. 
\item If $n=5$ then the foliation conjecture is true with weight vector $\mu_c=\left(\frac{1}{2}, \dots, \frac{1}{2} \right)$. 
\end{enumerate}
\end{thm}

We now proceed to describe briefly the contents of the paper. In Section~\ref{basic} we introduce our moduli spaces of parabolic vector bundles and Higgs bundles over the five punctured projective line, and  give some background on elementary transformations.  In Section~\ref{locusunst} we study the locus of Higgs fields which admit  unstable underlying parabolic vector bundle. Then, in Section~\ref{nilpotentcone}, we give an explicit description of the nilpotent cone, as well as the fixed points with respect to the $\C^*$-action. In Section~\ref{section:singnscurves},  we describe the remain singular fibers of the Hitchin map. Finally, in Section~\ref{section:connections} we introduce moduli spaces of connections and the foliation conjecture is investigated.

\section{Basic definitions}\label{basic}

Let $\Lambda=0+1+\lambda+t+\infty$ be a  divisor on the complex projective line $\mathbb P^{1}$ supported on five distinct points.  

\subsection{Moduli spaces}
A rank two \textit{quasiparabolic vector bundle} $ (E, {l})$, ${l} = \{l_{i}\}$,   on $\big(\P^1, \Lambda\big)$ consists  of a holomorphic vector bundle $E$ of rank two on $\P^1$ and for each  $i \in \{0, 1, \lambda, t, \infty\}$, a $1$-dimensional linear subspace $l_{i} \subset E_{i}$. We call $\Lambda$ the divisor of parabolic points, and the subspaces $l_{i} \subset E_i$ are called parabolic directions.

Let us now introduce a notion of stability for quasiparabolic vector bundles. Fix a weight  vector $\mu = (\mu_{1}, \dots, \mu_{5})$ of real numbers $0 \leq \mu_{i} \leq 1$.
A quasiparabolic vector bundle $(E,{l})$  is $\mu$-\textit{semistable} (respectively $\mu$-\textit{stable}) if for every  line subbundle $L \subset E$ we have 
\[
{\rm Stab}_{\mu}(L) := \deg E - 2\deg L - \sum_{l_i = L|_i} \mu_i + \sum_{l_i\neq L|_i} \mu_i \ge 0
\] 
(respectively the strict inequality holds). A \textit{parabolic vector bundle} is a quasiparabolic vector bundle together with a  weight vector $\mu$. We say that a parabolic vector bundle is {\it semistable} if the corresponding quasiparabolic vector bundle is $\mu$-\textit{semistable}. For each  $d\in\mathbb Z$ and a weight vector $\mu$, there is a {\it moduli space} $Bun_{\mu}(\P^1,\Lambda,d)$, 
parametrizing rank two parabolic vector bundles on $\big(\P^1, \Lambda\big)$, with $\deg E = d$,  which are semistable. 

Let us fix $d=0$. It follows from \cite{B}, that there is a polytope $\Delta\subset [0,1]^5$ consisting of weight vectors $\mu$ such that $Bun_{\mu}(\P^1,\Lambda,0)$ is nonempty. There are finitely many models $Bun_{\mu}(\P^1,\Lambda,0)$, corresponding to different chambers in the  wall-and-chamber decomposition of $\Delta$,  coming from to the variation  of the GIT. For example, the central weight $\mu_c = (\frac{1}{2},\dots, \frac{1}{2})$ lies in the interior of a chamber and  the moduli space 
\begin{eqnarray}\label{dels}
\mathcal S = Bun_{\mu_c}(\P^1, \Lambda, 0) 
\end{eqnarray}
is a del Pezzo surface of degree four, see also \cite{LS}.

 A {\it parabolic Higgs bundle} is a triple $(E,{l}, \theta)$ where $(E, {l})$ is a quasiparabolic vector bundle over $(\P^1, \Lambda)$ and  $\theta: E\to E\otimes\omega_{\P^1}(\Lambda)$ is a traceless homomorphism, which is nilpotent with respect to the parabolic directions. The condition of being nilpotent means that the residual part ${\rm Res}(\theta, i)$ satisfies ${\rm Res}(\theta, i)\cdot l_i = 0$ and ${\rm Res}(\theta, i)(E_i)\subset l_i$, for each $i\in\{0, 1, \lambda, t,\infty\}$.   We say that $\theta$ is a {\it parabolic Higgs field}.  A line subbundle $L\subset E$ is called {\it invariant} under $\theta$ if  $\theta(L)\subset L\otimes \omega_{\P^1}(\Lambda)$. In addition, $\theta$ is {\it irreducible} if it does not admit invariant line subbundle.


A parabolic Higgs bundle $(E, {l},\theta)$  is called $\mu$-\textit{semistable} (respectively $\mu$-\textit{stable}) if for every  line subbundle $L \subset E$ invariant under $\theta$, we have ${\rm Stab}_{\mu}(L) \ge 0$ (respectively ${\rm Stab}_{\mu}(L)>0$). 
We say that $(E, {l},\theta)$ is $\mu$-\textit{unstable} if it is not  $\mu$-\textit{semistable}. 

It follows from \cite[Propositions 3.1 and 3.2]{FL22}, that every parabolic Higgs field $\theta$ on $(\P^1, \Lambda)$ with $\det \theta \neq 0$ is irreducible, then $\mu$-stable for any choice of weight vector. Note also that the condition of being nilpotent implies that the quadratic differential $\det \theta$ lies in ${\rm H}^0(\P^1, \omega_{\P^1}^{\otimes 2}(\Lambda))$, which is a two dimensional vector space.  

For each weight vector $\mu$ there is a moduli space $\mathcal H_{\mu}(\P^1, \Lambda, 0)$ parametrizing  parabolic Higgs bundles on $(\P^1, \Lambda)$, with $\deg E = 0$, which are $\mu$-semistable \cite{Yo93, Yo95}. We denote by   
\begin{eqnarray}\label{sur P1}
\mathcal H = \mathcal H_{\mu_c} (\mathbb P^{1}, \Lambda,0)
\end{eqnarray}
the moduli space corresponding to the central weight $\mu_c$. It is a smooth four dimensional quasiprojective variety. 

%
%
%

\subsection{Elementary transformations}
The automorphism group of $\mathcal S$, cf. (\ref{dels}), has order $16$,  and admits a modular interpretation in terms of the group ${\bf El}$ of elementary transformations \cite{LS, AFKM21},  which we now describe.  

Assume that $I\subset \{0,1,\lambda, t, \infty\}$ has  even cardinality and let 
\[
D_I=\sum_{i\in I} i
\]
be the corresponding divisor. We consider the following exact sequence of sheaves
$$
0\ \rightarrow\ E' \ \stackrel{\alpha}{\to} \  E\ {\rightarrow}\ \bigoplus_{i\in I} E/l_i \ \rightarrow\ 0 \ 
$$
where $E/l_i$ intends to be a skyscraper sheaf determined by $E_{i}/l_i$.
We view $E'$ as a quasiparabolic vector bundle $(E',{l'})$ of rank two over $(\P^1, \Lambda)$ putting $l_i' := {\rm ker}\alpha_{i}$. We call it the {\it elementary transformation} of $(E,{l})$ over $D_I$:
\[
elem_{D_I}(E,{l}) := (E', {l'}).
\]
After this correspondence, the determinant line bundle is affected  
\[
\det E' = \det E \otimes \mathcal O_{\P^1}(-D_I),
\]
we then take  a square root $L_I$ of $\mathcal O_{\P^1}(D_I)$, in order to obtain
\[
\det (E'\otimes L_I) = \mathcal O_{\P^1}. 
\]

The stability condition is preserved if the weight vector $\mu=(\mu_1,\dots,\mu_5)$ is modified as follows.   If $(E,{l})$ is $\mu$-semistable then $elem_{D_I}(E,{l})$ is $\mu'$-semistable with $\mu_i'=\mu_i$ if $i\notin I$ and $\mu_i'=1-\mu_i$ if $i\in I$.  
In particular, when $\mu$ is the central weight, 
we obtain an isomorphism  
\begin{eqnarray}\label{eleI}
elem_{I}:\mathcal S \to \mathcal S
\end{eqnarray}
which sends $(E,{l})$ to $elem_{D_I}(E,{l})\otimes L_I$. 

It follows from basic properties of elementary transformations that 
\[
elem_{I} \circ elem_{J} = elem_{K}
\]
where $K=I\cup J\setminus I\cap J$ and the group ${\bf El}$ of transformations of the form  $elem_{I}$, where $I$ runs over all the subsets of $ \{0,1,\lambda, t, \infty\}$ of even cardinality, gives a modular realization of $\left(\frac{\mathbb Z}{2\mathbb Z}\right)^4$. Besides this, ${\bf El}$  coincides with the whole automorphism group of $\mathcal S$.

Note that, similarly, each correspondence $elem_I$  also acts on Higgs bundles, giving a modular realization of $\left(\frac{\mathbb Z}{2\mathbb Z}\right)^4$ as subgroup  of the automorphism group of $\mathcal H$, which we still denote by ${\bf El}\subset {\rm Aut}\mathcal H$.   See also \cite[Section 2.4]{FL22} and \cite[Section 4.2]{FL} for more details on elementary transformations.

\section{Higgs fields having unstable parabolic bundles}\label{locusunst}

Let $\mathcal S$ and $\mathcal H$ be as in the previous section.  There is an embedding $\mathcal  S\to \mathcal H$ by taking the Higgs field to be zero.  Since the weight vector $\mu_c=\left(\frac{1}{2}, \dots, \frac{1}{2} \right)$ lies in the interior of a chamber, any parabolic vector bundle in $\mathcal S$ is $\mu_c$-stable. It might happen that $(E, {l}, \theta)$ is $\mu$-semistable with  $(E, {l})$ $\mu$-unstable. For instance an unstable parabolic bundle may be often endowed with an irreducible Higgs field.  In this  section we shall studying this phenomenon.

Let us consider the forgetful map 
\[
\frak{for} : \mathcal H \dashrightarrow \mathcal S
\]
which forgets the Higgs field.  There is an open subset of $\mathcal H$ where $\frak{for}$ is well defined, it is formed by Higgs bundles over $\mathcal S$: 
\[
\mathcal U = \{ (E, {l}, \theta)\in \mathcal H\;:\;\; (E, {l})\in \mathcal S\}.
\]
There is an identification between $\mathcal U$ and the cotangent bundle $T^*\mathcal S$, by identifying $T^*_{(E, {l})}\mathcal S$ with $\frak{for}^{-1}(E, {l})$, see \cite[Theorem 2.4]{Yo95}, so $\mathcal H$  contains the cotangent bundle $T^*\mathcal S$ as an open and dense subset.


The next result describes which underlying parabolic bundles appear in $\mathcal H$. 
 
\begin{prop}\label{classunder}
Given $(E, {l}, \theta)\in \mathcal H$, then 
\begin{itemize}
\item $E = \mathcal O_{\P^1}(-d)\oplus  \mathcal O_{\P^1}(d)$, with $d\in \{0,1\}$;
\item if $d=0$ then at most $3$ parabolic directions lie in the same embedding of $\mathcal O_{\P^1}\hookrightarrow E$;
\item if $d=1$ then at most $1$ parabolic direction lies in  $\mathcal O_{\P^1}(1)\hookrightarrow E$. 
\end{itemize}
\end{prop}
 
\proof
Since $E$ has degree zero, we can assume $E   = \mathcal O_{\P^1}(-d)\oplus  \mathcal O_{\P^1}(d)$, with $d\ge 0$. A Higgs field 
\begin{eqnarray*}
\theta=\left(
\begin{array}{ccc} 
\alpha & \beta  \\
\gamma & -\alpha  \\
\end{array}
\right)
\end{eqnarray*}
with logarithmic poles on $\Lambda$ is given by homomorphisms 
\begin{displaymath}
\left\{ \begin{array}{ll}
\alpha: \mathcal O_{\P^1} \to \omega_{\P^1}(\Lambda)\\
\beta: \mathcal O_{\P^1}(d) \to \mathcal O_{\P^1}(-d)\otimes\omega_{\P^1}(\Lambda)\\
\gamma: \mathcal O_{\P^1}(-d) \to \mathcal O_{\P^1}(d)\otimes\omega_{\P^1}(\Lambda)
\end{array} \right.
\end{displaymath}
which turns out to be equivalent to give
\begin{displaymath}
\left\{ \begin{array}{ll}
\alpha \in \Gamma(\mathcal O_{\P^1}(3))\\
\beta \in \Gamma(\mathcal O_{\P^1}(3-2d))\\
\gamma \in \Gamma(\mathcal O_{\P^1}(3+2d))
\end{array} \right. 
\end{displaymath}
Now if $d\ge 2$ then $\beta=0$ and  $\mathcal O_{\P^1}(d)$ is a destabilizing subbundle. This concludes the first assertion of the statement. 

Let us assume $d=0$. An embedding  $\mathcal O_{\P^1}\hookrightarrow \mathcal O_{\P^1}\oplus\mathcal O_{\P^1}$, $e\mapsto (e,0)$,   passing through a parabolic direction $l_i$ over $t_i$ yields $\gamma\in\Gamma(\mathcal O_{\P^1}(3)\otimes\mathcal O_{\P^1}(-t_i))$. Thus at most $3$ parabolic directions lie in $\mathcal O_{\P^1}$, otherwise $\gamma=0$ and $\mathcal O_{\P^1}$ is a destabilizing subbundle. The case $d=1$ is similar, and hence will be  omitted.    
\endproof 

 \begin{cor}\label{cor:under}
Let  $(E, {l}, \theta)\in \mathcal H$. Assume that the underlying parabolic bundle $(E, {l})$ is $\mu_c$-unstable. Then we are in one of the following possibilities
\begin{itemize}
\item $E=L_1\oplus L_2$, $L_i\simeq \mathcal O_{\P^1}$, $L_1$ contains $3$ parabolic directions and $L_2$ contains $2$ parabolic directions; 
\item $E=\mathcal O_{\P^1}(-1)\oplus \mathcal O_{\P^1}(1)$ and $\mathcal O_{\P^1}(-1)$ contains every parabolic direction;
\item $E=\mathcal O_{\P^1}(-1)\oplus \mathcal O_{\P^1}(1)$,  $\mathcal O_{\P^1}(1)$ contains exactly $1$ parabolic direction and $\mathcal O_{\P^1}(-1)$ contains the remaining $4$ parabolic directions. 
\end{itemize}
In particular,  $(E, {l})$ is decomposable. 
 \end{cor}
 
 \proof
We first  reduce to the case where $E$ is trivial, up to an elementary transformation. If $E=\mathcal O_{\P^1}(-1)\oplus \mathcal O_{\P^1}(1)$, Proposition~\ref{classunder} ensures that at most $1$ parabolic directions lies in $\mathcal O_{\P^1}(1)$, and since the family of embeddings $\mathcal O_{\P^1}(-1)\hookrightarrow E$ is three dimensional, we can take an  $\mathcal O_{\P^1}(-1)$ passing through $3$ parabolic directions outside $\mathcal O_{\P^1}(1)$. Now,  a transformation $elem_I$ over two of them, transforms $E$ into the trivial vector bundle. 

Assume that $E$ is trivial and $(E, {l})$ is $\mu_c$-unstable. A destabilizing subbundle $L\subset E$, $\deg L\le 0$, satisfies 
\[
-2\deg L-m/2+n/2<0
\]
where $m$ is the number of parabolic directions in $L$ and $n$ corresponds to the parabolic directions outside $L$. Hence, $\deg L\in \{0,-1\}$. In addition, by Proposition~\ref{classunder}, if $\deg L=0$ then there are exactly $3$ parabolic directions in $L$.  If $\deg L = -1$ then every parabolic direction lies in $L$, and applying  a transformation $elem_I$ over two parabolic points,  we reduce to the previous case. 

Now we may assume that $E$ is trivial, and there are exactly $3$  parabolic directions in the same embedding $\mathcal O_{\P^1}\hookrightarrow L_1\subset E$. We will show that $\mu_c$-semistability of $\theta$ implies  that there exists another embedding of $\mathcal O_{\P^1}\hookrightarrow L_2\subset E$ passing through the remaining two parabolic directions. 

By simplicity, let us assume that parabolic directions $l_0, l_1, l_{\lambda}$ over $0, 1, \lambda$ lie in $L_1$ and let $L_2$ be an embedding of $\mathcal O_{\P^1}$ passing through the parabolic direction $l_t$ and we have $E=L_1\oplus L_2$. As in the proof of Proposition~\ref{classunder}, since $\theta$ is nilpotent with respect to the parabolic directions, then $\gamma$ vanishes at $\{0, 1, \lambda\}$, $\beta$ vanishes at $t$, and $\alpha$ vanishes at $\{0, 1, \lambda, t\}$. So, we conclude that $\alpha=0$ and 
\begin{displaymath}
\left\{ \begin{array}{ll}
\beta: \mathcal O_{\P^1} \to \omega_{\P^1}(0+1+\lambda+\infty)\\
\gamma: \mathcal O_{\P^1}\to \omega_{\P^1}(t+\infty).
\end{array} \right.
\end{displaymath}
If the remaining parabolic direction $l_{\infty}$ over $\infty$ is outside $L_2$ then the condition to be nilpotent implies that $\beta$ and $\gamma$ vanish  on it. In this case, $\gamma$ must be zero, $L_1$ is invariant under  $\theta$ and then $\theta$ is $\mu_c$-unstable. When it lies in $L_2$ then $\beta$ vanishes also at $\infty$, i.e., $\beta: \mathcal O_{\P^1} \to \omega_{\P^1}(0+1+\lambda)$.

We have shown that $E=L_1\oplus L_2$, $L_i=\mathcal O_{\P^1}$,  $3$ parabolic directions $l_0, l_1, l_{\lambda}$ lie in $L_1$,  and the remaining directions $l_t, l_{\infty}$ lie in $L_2$. This concludes the proof of the corollary. 
 
 \endproof

\begin{center}
\begin{figure}[h]
\centering
\includegraphics[height=2.2in]{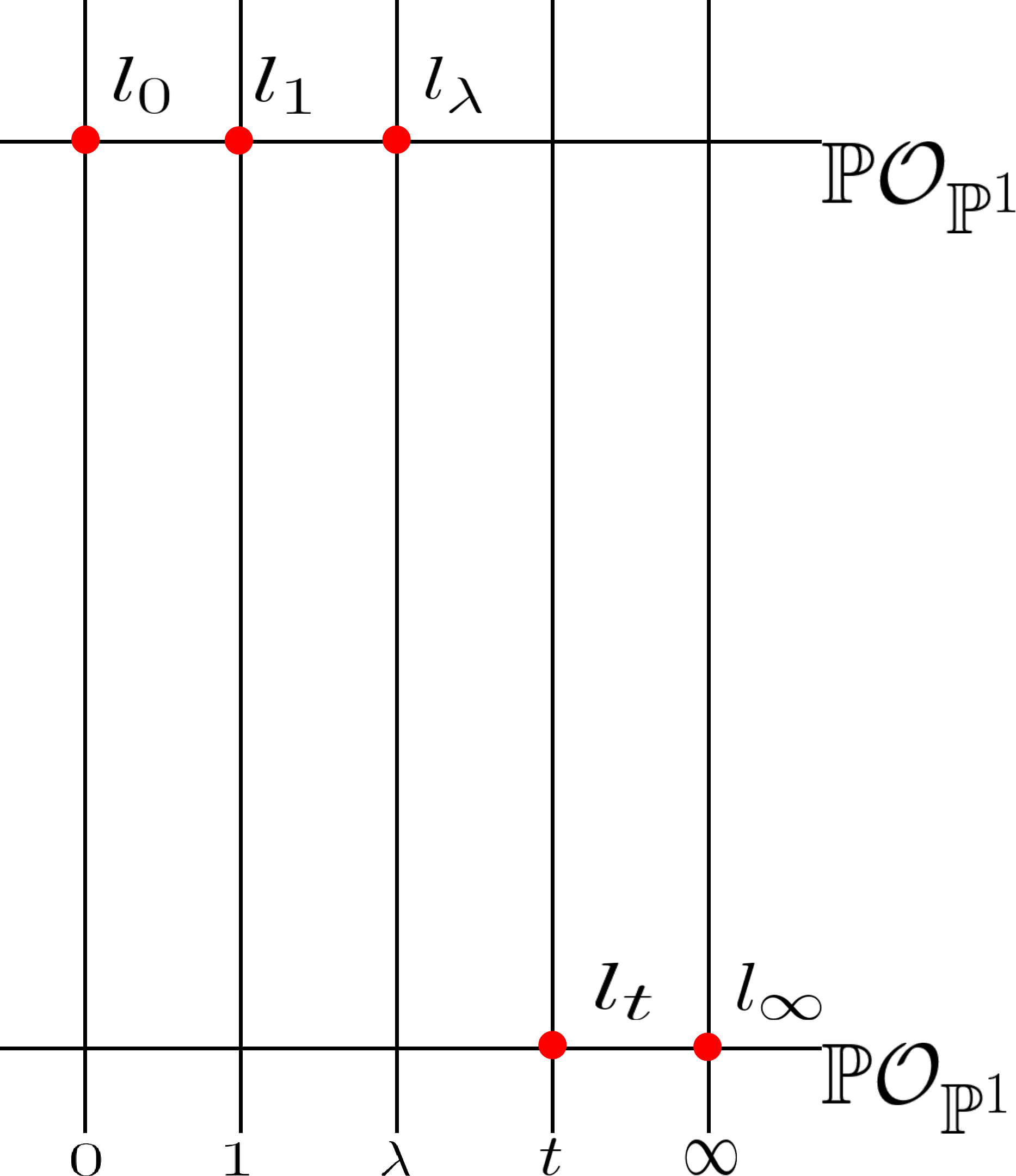}
\caption{Unstable parabolic bundle which admits stable Higgs field.}
\label{parunst}
\end{figure}
\end{center}

This corollary implies that there are exactly $16$ $\mu_c$-unstable parabolic vector bundles which admit a $\mu_c$-semistable Higgs field $\theta$, see Table~\ref{16}.  The group ${\bf El}$ acts transitively on it and  Figure~\ref{parunst} shows one of them. 

  \begin{table}
\centering
\begin{tabular}[p]{|c|c|c|}
   \hline
            & $E$ &  $\{u,v,p,q,r\}=\{0, 1, \lambda, t, \infty\}$\\
   \hline \hline           
   $10$ & $\mathcal O_{\P^1}\oplus\mathcal O_{\P^1}$ & $l_u, l_v, l_p\subset L_1\simeq \mathcal O_{\P^1}$ and $l_q, l_r\subset L_2\simeq \mathcal O_{\P^1}$\\
   \hline
   $5$ & $\mathcal O_{\P^1}(-1)\oplus \mathcal O_{\P^1}(1)$ & $l_u\subset \mathcal O_{\P^1}(1)$ and $l_v, l_p, l_q, l_r\subset \mathcal O_{\P^1}(-1)$ \\
   \hline
   $1$ & $\mathcal O_{\P^1}(-1)\oplus \mathcal O_{\P^1}(1)$ & $l_0, l_1, l_{\lambda}, l_t, l_{\infty}\subset \mathcal O_{\P^1}(-1)$\\
   \hline
  \end{tabular}
\caption{$16$ unstable parabolic bundles admitting  stable Higgs fields}
\label{16}
\end{table}

\begin{remark}\rm
Let us assume we are in the first case of Corollary~\ref{cor:under}:  $E=L_1\oplus L_2$, $L_i\simeq \mathcal O_{\P^1}$, $L_1$ contains $3$ parabolic directions, over $0, 1$ and $\lambda$, and $L_2$ contains $2$ parabolic directions, over $\infty$ and  $t$. We have seen that any $\mu_c$-semistable  Higgs field on it is of the form  
\begin{eqnarray*}
\theta=\left(
\begin{array}{ccc} 
0 & \beta  \\
\gamma & 0 \\
\end{array}
\right)
\end{eqnarray*}
with 
\begin{displaymath}
\left\{ \begin{array}{ll}
\beta: \mathcal O_{\P^1} \to \omega_{\P^1}(0+1+\lambda)\\
\gamma: \mathcal O_{\P^1}\to \omega_{\P^1}(t+\infty)\;, \quad \gamma\neq 0\,.
\end{array} \right.
\end{displaymath}
Any other Higgs bundle admitting a $\mu_c$-unstable parabolic vector bundle can be obtained from this by performing an elementary transformation.  
\end{remark}

In the next result we will determine the complement $\mathcal H\setminus \mathcal U$,   formed by $\mu_c$-semistable Higgs bundles which have $\mu_c$-unstable underlying parabolic bundle.  Before that, let us introduce some notation: let ${\rm Higgs}(E,{l})$ be the quotient of the vector space $\Gamma(\mathcal{SE}nd( E, {l})\otimes \omega_{\P^1}(\Lambda))$ by the automorphism group of the parabolic bundle $(E,{l})$. The stability condition has not been considered here, a point of ${\rm Higgs}(E,{l})$ lies in $\mathcal H$ only if it is $\mu_c$-semistable. 



\begin{prop}\label{prop:16irred}
The complement $\mathcal H\setminus \mathcal U$  has exactly $16$ irreducible components and the group ${\bf El}$ acts transitively on it. Each component is a Zariski open subset of ${\rm Higgs}(E,{l})$, for each one of the $16$ decomposable parabolic bundles shown in Table~\ref{16}. 
\end{prop}

\proof
An element $(E, {l}, \theta)$ of $\mathcal H\setminus \mathcal U$ corresponds to a Higgs field which has $\mu_c$-unstable underlying parabolic bundle. These parabolic bundles were classified in Corollary~\ref{cor:under} and there are $16$ of them. In addition, the group ${\bf El}$ acts transitively on it, so we fix one, say  $E=L_1\oplus L_2$, $L_i=\mathcal O_{\P^1}$,  with $3$ parabolic directions over $0, 1, \lambda$ lying  in $L_1$,  and with the remaining directions, over $t, \infty$,  lying in $L_2$. The corresponding space of Higgs fields 
\[
 \Gamma(\mathcal{SE}nd( E, {l})\otimes \omega_{\P^1}(\Lambda))
\]
is three dimensional and its quotient by the automorphism group of $( E, {l})$ gives  ${\rm Higgs}(E,{l})$. We want the locus in ${\rm Higgs}(E,{l})$ formed by $\mu_c$-semistable Higgs fields. According to the proof of Corollary~\ref{cor:under},  any  Higgs field in ${\rm Higgs}(E,{l})$ is given by
\begin{eqnarray}\label{betagamma}
\theta=\left(
\begin{array}{ccc} 
0 & \beta  \\
\gamma & 0 \\
\end{array}
\right)
\end{eqnarray}
where 
\begin{displaymath}
\left\{ \begin{array}{ll}
\beta: \mathcal O_{\P^1} \to \omega_{\P^1}(0+1+\lambda)\\
\gamma: \mathcal O_{\P^1}\to \omega_{\P^1}(t+\infty)
\end{array} \right.
\end{displaymath}
and so $(\beta, \gamma)$ lies in a three dimensional vector space. We see that $\theta$  is $\mu_c$-semistable if and only if $\gamma\neq 0$.  On the other hand, automorphisms of $(E,{l})$, i.e. automorphisms of the trivial bundle fixing parabolic directions, are diagonal and then the quotient of 
\[
\Gamma(\mathcal{SE}nd( E, {l})\otimes \omega_{\P^1}(\Lambda))\setminus \{\gamma=0\}
\]
is a  two dimensional subvariety of $\mathcal H$.  
\endproof

\section{Nilpotent cone}\label{nilpotentcone}

The nilpotent cone $\mathcal N$ is formed by Higgs fields having vanishing determinant, we will show that it has $17$ irreducible components. Of course it contains $\mathcal S$, the locus obtained by taking $\theta=0$, which is a del Pezzo surface of degree four. We will show that outside $\mathcal S$  there is exactly one component for each of the $16$ special rational curves of $\mathcal S$ (those which have $(-1)$-self intersection).  These curves are parametrized by  parabolic structures given in Table~\ref{16del}, see \cite{LS} for details.

 \begin{table}
\centering
\begin{tabular}[p]{|c|c|c|}
   \hline
            & $E$ &  $\{u,v,p,q,r\}=\{0, 1, \lambda, t, \infty\}$\\
   \hline \hline           
   $10$ & $\mathcal O_{\P^1}\oplus\mathcal O_{\P^1}$ & $l_u, l_v\subset \mathcal O_{\P^1}\hookrightarrow E$ \\
   \hline
   $5$ & $\mathcal O_{\P^1}\oplus \mathcal O_{\P^1}$ & $l_u, l_v, l_p, l_q\subset \mathcal O_{\P^1}(-1)\hookrightarrow E$ \\
   \hline
   $1$ & $\mathcal O_{\P^1}(-1)\oplus \mathcal O_{\P^1}(1)$ & $l_0, l_1, l_{\lambda}, l_t, l_{\infty}\nsubseteq \mathcal O_{\P^1}(1)$\\
   \hline
  \end{tabular}
\caption{$16$ special lines in $\mathcal S$}
\label{16del}
\end{table}

We first determine the intersection between $\mathcal N$ and $ \mathcal H\setminus \mathcal U$, i.e., $\mu_c$-semistable Higgs bundles having $\mu_c$-unstable parabolic vector bundle and with vanishing determinant. To give one example, let 
\[
\Theta_1= (L_1\oplus L_2, {l}, \theta)\;,\quad L_i\simeq \mathcal O_{\P^1}
\] 
where the parabolic structure is given by 
\[
l_0, l_1, l_{\lambda}\subset  L_1\quad\text{and}\quad  l_t, l_{\infty}\subset L_2
\] 
and
\begin{eqnarray}\label{ex:theta0}
\theta=\left(
\begin{array}{ccc} 
0 & 0 \\
\frac{dx}{(x-t)}  & 0 \\
\end{array}
\right).
\end{eqnarray}
Note that the destabilizing subbundle $L_1$ for the underlying parabolic structure is non-invariant under $\theta$. By performing the transformations $elem_I\in {\bf El}$  we get at least $16$ $\mu_c$-semistable Higgs bundles, $\Theta_i$, $i=1,\dots, 16$,  having $\mu_c$-unstable underlying parabolic vector bundles. The next result shows that these are all the cases.

\begin{prop}\label{prop:Theta_i}
There are exactly $16$ Higgs bundles in $ \mathcal H\setminus \mathcal U$ with vanishing determinant,  they are $\Theta_i$, $i=1,\dots, 16$, as above. 
%
\end{prop}  

\proof

By Proposition~\ref{prop:16irred} we may assume, up to a transformation $elem_I$,  that a Higgs bundle $(E, {l}, \theta)\in\mathcal H\setminus \mathcal U$ is given by $\theta$ as in (\ref{betagamma}) 
and $(E, {l})$ is the parabolic vector bundle of Figure~\ref{parunst}. Now, if $\theta$ has  vanishing determinant then $\beta\gamma=0$, and $\gamma$ cannot be zero because otherwise $\theta$ is $\mu_c$-unstable. We conclude that $\beta=0$ and  up to an automorphism of $(E, {l})$ we can assume that $\gamma$ has residue $1$ at $t$. This gives the expression for $\theta$ in (\ref{ex:theta0}). 

\endproof 

Let us denote by $\zeta_i\subset \mathcal S$, $i=1, \dots, 16$, the $(-1)$-self intersection  rational curves in $\mathcal S$, see Table~\ref{16del}, and let 
\[
\Sigma = \cup_{i=1}^{16} \zeta_i
\] 
be the union. There is a natural correspondence between the set of rational curves
\[
\{\zeta_i\;:\; i=1,\dots, 16\}
\]
and the set of Higgs bundles 
\[
\{ \Theta_i\;:\; i=1,\dots, 16\}
\]
in $\mathcal N\cap (\mathcal H\setminus \mathcal U)$.
For instance, we first associate the rational curve $\zeta_1\subset \mathcal S$, corresponding to  parabolic vector bundles with two parabolic directions $l_t, l_{\infty}$ inside $\mathcal  O_{\P^1}\simeq L_2\subset \mathcal O_{\P^1}\oplus\mathcal O_{\P^1}$, to $\Theta_1$. The structure of the underlying parabolic vector bundle of $\Theta_1$ is infinitely close to the parabolic structure varying in $\zeta_1$.
%
Now, the correspondence 
\[
\zeta_i \longleftrightarrow \Theta_i
\]
follows by the action of ${\bf El}$ in both sets. 

We will see that besides $\mathcal S$, the nilpotent cone has $16$ components $\mathcal N_i$ which can be obtained as one point compactification of $\mathcal N_i\cap \mathcal U$, i.e.
\[
\mathcal N_i = (\mathcal N_i\cap \mathcal U)\cup \{\Theta_i\}
\]
So, first we study the restriction of the nilpotent cone to $\mathcal U$. 

\begin{prop}\label{prop:sigmaimage}
If $(E, {l}, \theta)\in \mathcal U$ has vanishing determinant, then $(E, {l})\in \Sigma$, ie. $\frak{for}(\mathcal N) = \Sigma$.  
\end{prop}
 
\proof
To begin with, note that if $E=\mathcal O_{\P^1}(-1)\oplus \mathcal O_{\P^1}(1)$,  we can apply  $elem_I$ in order to transform $E$ into the trivial vector bundle $E=\mathcal O_{\P^1}\oplus \mathcal O_{\P^1}$. In addition, since $(E, {l})$ is $\mu_c$-semistable, there is no embedding of $\mathcal O_{\P^1}\hookrightarrow E$  passing through $3$ parabolic directions, and then  for computation we can assume that the parabolic directions ${l} = \{l_i\}$ are normalized as 
\[
l_0= \begin{pmatrix}1\\ 0\end{pmatrix}, l_1= \begin{pmatrix}1\\ 1\end{pmatrix}, l_{\lambda}= \begin{pmatrix}1\\ u\end{pmatrix}, l_t= \begin{pmatrix}1\\ v\end{pmatrix},  l_{\infty}= \begin{pmatrix}0\\ 1\end{pmatrix}.
\]

Any Higgs field $\theta$ on $(E, {l})$ can be written as 
\[
\theta  = c_1\theta_1 + c_2\theta_2 \;;  \quad c_1, c_2\in \mathbb C
\]
where 
\begin{eqnarray*}
\theta_1=\left(
\begin{array}{ccc} 
\frac{u}{(x-\lambda)} - \frac{u}{(x-1)}   & \quad \frac{u}{(x-1)} - \frac{1}{(x-\lambda)} +\frac{1-u}{x} \\
\frac{u^2}{(x-\lambda)} - \frac{u}{(x-1)} & \quad  -\frac{u}{(x-\lambda)} + \frac{u}{(x-1)} \\
\end{array}
\right)\cdot dx
\end{eqnarray*}
and 
\begin{eqnarray*}
\theta_2=\left(
\begin{array}{ccc} 
\frac{v}{(x-t)} - \frac{v}{(x-1)}   & \quad \frac{v}{(x-1)} - \frac{1}{(x-t)} +\frac{1-v}{x} \\
\frac{v^2}{(x-t)} - \frac{v}{(x-1)} & \quad  -\frac{v}{(x-t)} + \frac{v}{(x-1)} \\
\end{array}
\right)\cdot dx
\end{eqnarray*}
and $x$ intends to be the coordinate of $\P^1$. Then we get 
\[
\det \theta = (h_1+h_2\cdot x)\frac{dx^{\otimes 2}}{x(x-1)(x-\lambda)(x-t)}
\]
where
\begin{eqnarray*}
h_1 &=& (c_1(1-u)+c_2(1-v))(c_1tu (\lambda -u) + c_2\lambda v(t-v))  \\
h_2 &=& (c_1u(u-1)+c_2v(v-1))(c_1(\lambda-1)+c_2(t-v))
\end{eqnarray*}
and  let us write 
\begin{displaymath}
\left\{ \begin{array}{ll}
a_{1} = c_1(1-u)+c_2(1-v)\\
a_{2}=c_1tu (\lambda -u) + c_2\lambda v(t-v)\\
b_{1}=c_1u(u-1)+c_2v(v-1)\\
b_{2}=c_1(\lambda-1)+c_2(t-v)
\end{array} \right.
\end{displaymath}
We see that $\theta$ has vanishing determinant if and only if 
\begin{eqnarray}\label{linsystem}
a_i = b_j =0 
\end{eqnarray}
for some $i,j\in \{1,2\}$. 

We are looking for nontrivial solutions  $c_1,c_2$ for each linear system (\ref{linsystem}) and actually we will show that  it has a nontrivial solution if and only if the parabolic structure lies in $\Sigma$, which is  the locus of $16$ special rational curves of $\mathcal S$. To do so, we first note that the system  $a_i=b_j=0$, for some $i,j\in\{1,2\}$,   has a nontrivial solution $c_1,c_2$ if and only if at least one of the following equations hold
\begin{displaymath}
\left\{ \begin{array}{ll}
 (v-1)(u-1)(u-v) = 0\\
(t-v)(-u+\lambda)(\lambda v-t u) = 0\\
u(t-1)+v(1-\lambda) +\lambda-t = 0\\
vu( u t (\lambda-1)+v \lambda (1-t)+uv(t-\lambda)) = 0
\end{array} \right.
\end{displaymath}
This last  means that either there are two parabolic points lying in the same embedding $\mathcal O_{\P^1}\hookrightarrow \mathcal O_{\P^1}\oplus \mathcal O_{\P^1}$ or there is an embedding $\mathcal O_{\P^1}(-1)\hookrightarrow \mathcal O_{\P^1}\oplus \mathcal O_{\P^1}$ passing through $4$ parabolic directions. More precisely, there is an embedding $\mathcal O_{\P^1}(-1)\hookrightarrow \mathcal O_{\P^1}\oplus \mathcal O_{\P^1}$ passing through $l_0, l_1, l_{\lambda}, l_{\infty}$ when
\[
\lambda-u = 0,
\]
through $l_0, l_1, l_{t}, l_{\infty}$ when
\[
t-v=0,
\]
through $l_0, l_{\lambda}, l_t, l_{\infty}$ when
\[
\lambda v-t u=0,
\]
through $l_1, l_{\lambda}, l_t, l_{\infty}$ when
\[
u(t-1)+v(1-\lambda) +\lambda-t = 0,
\]
and through $l_0, l_1, l_{\lambda}, l_t$ when 
\[
 u t (\lambda-1)+v \lambda (1-t)+uv(t-\lambda) = 0.
\]
The other cases are evident. This shows that $(E, {l})\in \Sigma$, completing the proof of the proposition. 

\endproof 
 


We now emphasize another consequence of this proposition. Let  $\mathcal N_{(E,{l})}$ be the set of Higgs fields having $(E, {l})$ as underlying  parabolic bundle and having vanishing determinant. Along the proof of Proposition~\ref{prop:sigmaimage}, we have seen that the intersection $\mathcal N_{(E,{l})}\cap \mathcal U$ corresponds to a union of lines in the vector space 
\[
{\rm Higgs}(E, {l})\simeq \C^2.
\]  
More precisely, it is one single line when $(E, {l})\in \zeta_i$ and  $(E, {l})\notin \zeta_j$ for $j\neq i$, and exactly two lines when $(E, {l})\in \zeta_i\cap \zeta_j$.

For convenience we give an explicit example. All the other cases can be obtained from this with a transformation $elem_I$.  If  $(E, {l})\in \zeta_1$, i.e., $E=\mathcal O_{\P^1}\oplus \mathcal O_{\P^1}$ and
\[
l_0= \begin{pmatrix}1\\ 0\end{pmatrix}, l_1= \begin{pmatrix}1\\ 1\end{pmatrix}, l_{\lambda}= \begin{pmatrix}1\\ u\end{pmatrix}, l_t=   l_{\infty}= \begin{pmatrix}0\\ 1\end{pmatrix}
\]
then 
\[
\mathcal N_{(E,{l})}\cap \mathcal U = \{ (E, {l}, c\cdot\theta_1)\; : \quad c\in \C\}
\]
where 
\begin{eqnarray}\label{ex:theta1}
\theta_1=\left(
\begin{array}{ccc} 
0 & 0 \\
\frac{dx}{(x-t)}  & 0 \\
\end{array}
\right)
\end{eqnarray}
when $(E, {l})\notin \zeta_j$ for $j\neq 1$. But if  $(E, {l})$ lies in the intersection of two rational curves $\zeta_1\cap \zeta_j$, for instance if $u=0$, then 
\[
\mathcal N_{(E,{l})}\cap \mathcal U = \{ (E, {l}, c\cdot\theta_1)\; : \quad c\in \C\}\cup  \{ (E, {l}, c\cdot\theta_j)\; : \quad c\in \C\}
\]
where
\begin{eqnarray}\label{ex:thetaj}
\theta_j=\left(
\begin{array}{ccc} 
0 & \frac{dx}{x(x-\lambda)}  \\
0 & 0 \\
\end{array}
\right).
\end{eqnarray}


It is interesting to note that for every $(E, {l})\in \zeta_1$ the line 
\[
 \{ (E, {l}, c\cdot\theta_1)\; : \quad c\in \C\}\subset \mathcal U
 \]
has the same limit point in $\mathcal H\setminus \mathcal U$, that is, 
\[
\lim_{c\to \infty} (E, {l}, c\cdot\theta_1) = \Theta_1
\]
where
\[
\Theta_1= (L_1\oplus L_2, {l}, \theta_1)\;,\quad L_i\simeq \mathcal O_{\P^1}
\] 
has  parabolic structure 
\[
l_0, l_1, l_{\lambda}\subset  L_1\quad\text{and}\quad  l_t, l_{\infty}\subset L_2.
\] 
In fact, for any $c\neq 0$, by performing an automorphism  
\begin{eqnarray}
\phi_c=\left(
\begin{array}{ccc} 
1 & 0 \\
0  & c^{-1} \\
\end{array}
\right)
\end{eqnarray}
on $(E, {l})$, one obtains 
\[
l_0= \begin{pmatrix}1\\ 0\end{pmatrix}, l_1= \begin{pmatrix}1\\ c^{-1}\end{pmatrix}, l_{\lambda}= \begin{pmatrix}1\\ c^{-1}u\end{pmatrix}, l_t=   l_{\infty}= \begin{pmatrix}0\\ 1\end{pmatrix}
\]
as parabolic directions, and hence when $c\to \infty$ the parabolic structure goes to the parabolic structure of $\Theta_1$. On the other hand, we have
\[
\phi\circ (c\cdot\theta_1)\circ\phi^{-1} = \theta_1.
\] 

Let $\theta_j$, $j=1,\dots, 16$, denote the transformed of $\theta_1$ by action of ${\bf El}$. We summarise the discussion above in the next result.

\begin{prop}\label{prop:lim}
If $(E, {l})\in \Sigma$ belongs to the rational curve $\zeta_i$, then 
\[
\mathcal N_{(E,{l})}\cap \mathcal U = \{ (E, {l}, c\cdot \theta_i)\; : \quad c\in \C\} 
\]
when $(E, {l})\notin \zeta_j$, $\forall j\neq i$,  and 
\[
\mathcal N_{(E,{l})}\cap \mathcal U = \{ (E, {l}, c\cdot \theta_i)\; : \quad c\in \C\} \cup \{ (E, {l}, c\cdot \theta_j)\; : \quad c\in \C\} 
\]
when $(E, {l})\in \zeta_i\cap \zeta_j$. Moreover, we have
\[
\lim_{c\to \infty} (E, {l}, c\cdot\theta_i) = \Theta_i\;.
\]
\end{prop}

\begin{center}
\begin{figure}[h]
\centering
\includegraphics[height=3.1in]{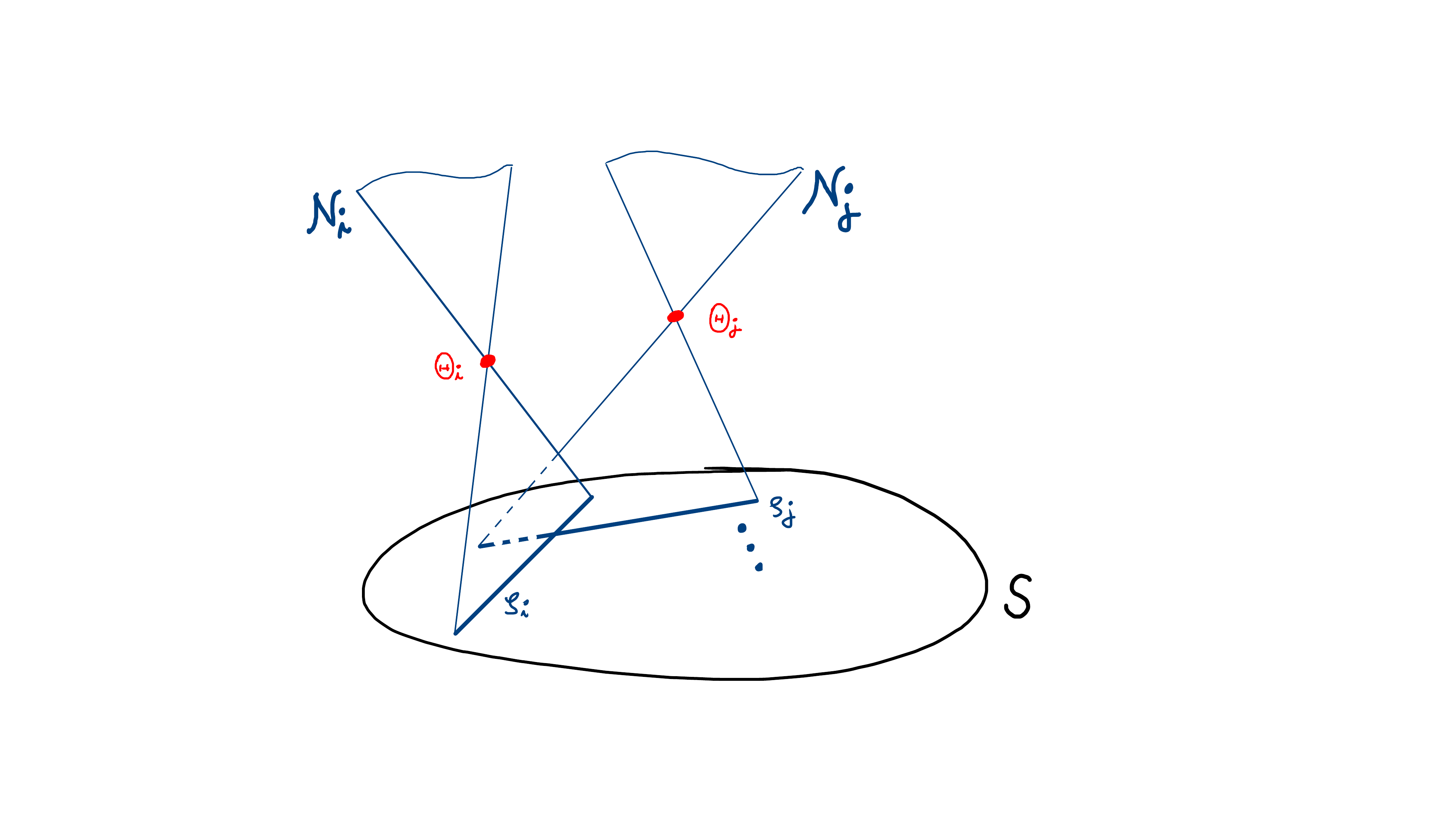}
\caption{Nilpotent cone.}
\label{nilpcone}
\end{figure}
\end{center}

 
\begin{definition}\label{def:Hodge}
The $\Theta_i$ are fixed points by the $\C^*$ action,  following the terminology of \cite{S08}, we call them the $16$ {\it Hodge bundles} of $\mathcal H$.  They are all the fixed points outside $\mathcal S$. 
\end{definition}
Finally, we are ready to the main result of this section:   
 
\begin{thm}\label{thm:nilpP1}
The nilpotent cone of $\mathcal H$ has exactly $17$ irreducible components
\[
\mathcal N = \mathcal S\cup_{i=1}^{16}\mathcal N_i
\]
where 
\[
\mathcal N_i = \{ (E, {l}, c\cdot \theta_i)\; : \;(E, {l})\in \zeta_i ,\;c\in \C\}\cup \{\Theta_i\}.
\]
See Figure~\ref{nilpcone}.
\end{thm} 
 
\proof
The proof follows from Propositions~\ref{prop:Theta_i}, \ref{prop:sigmaimage} and \ref{prop:lim}. 
\endproof

\section{Other singular Hitchin fibers}\label{section:singnscurves}

In this section we study singular fibers of the Hitchin map
\[
\det:\mathcal H \to \Gamma(\omega_{\P^1}^{\otimes 2}(\Lambda))\simeq \C^2
\]
over a point $s\neq 0$. The general spectral curve $X_s$ is a smooth curve of genus $2$ branched over $6$ distinct points 
\[
0,1,\lambda, t, \infty, \rho
\]
 of $\P^1$ and the corresponding Hitchin fiber is  ${\rm Pic}^{3}(X_s)$. A singular spectral curve occurs when the sixth point $\rho$ coincides with one of the five other points. Hence, the {\it locus of singular spectral curves} is a union of five lines
\[
\cup_{\rho} \Gamma(\omega_{\P^1}^{\otimes 2}(\Lambda-\rho)) 
\]
where $\rho$ varies in $\{0,1,\lambda, t, \infty\}$. If $s\neq 0$ lies in one of these lines then $X_s$ is a nodal curve of genus $2$, its desingularization $\tilde{X_s}$ is an elliptic curve branched over 
\[
\{0,1,\lambda, t, \infty\}\setminus\{\rho\} 
\]
and $X_s$ can be obtained identifying two points $w_{\rho}^+$ and $w_{\rho}^-$ of $\tilde{X_s}$. 

\begin{remark}\rm\label{rmk:Jac}
When $X_s$ is a nodal curve with a single node at $w_{\rho}$, its compactified Jacobian $\overline{{\rm Pic}}^0(X_s)$ is obtained identifying the $0$-section with the $\infty$-section, see Figure~\ref{p1bundle}, of the $\P^1$-bundle 
\begin{eqnarray}\label{compJac}
{\bf F} = \P(\mathcal O_{\tilde{X_s}}(w_{\rho}^+)\oplus\mathcal O_{\tilde{X_s}}(w_{\rho}^-))
\end{eqnarray}
via the translation $\mathcal O_{\tilde{X_s}}(w_{\rho}^+ - w_{\rho}^-)$, see (cf. \cite[p. 83]{OS79}). In particular, we have
\[
\tilde{X_s}\simeq \overline{{\rm Pic}}^0(X_s)\setminus {\rm Pic}^0(X_s).
\]
\hfill$\lrcorner$
\begin{center}
\begin{figure}[h]
\centering
\includegraphics[height=1.7in]{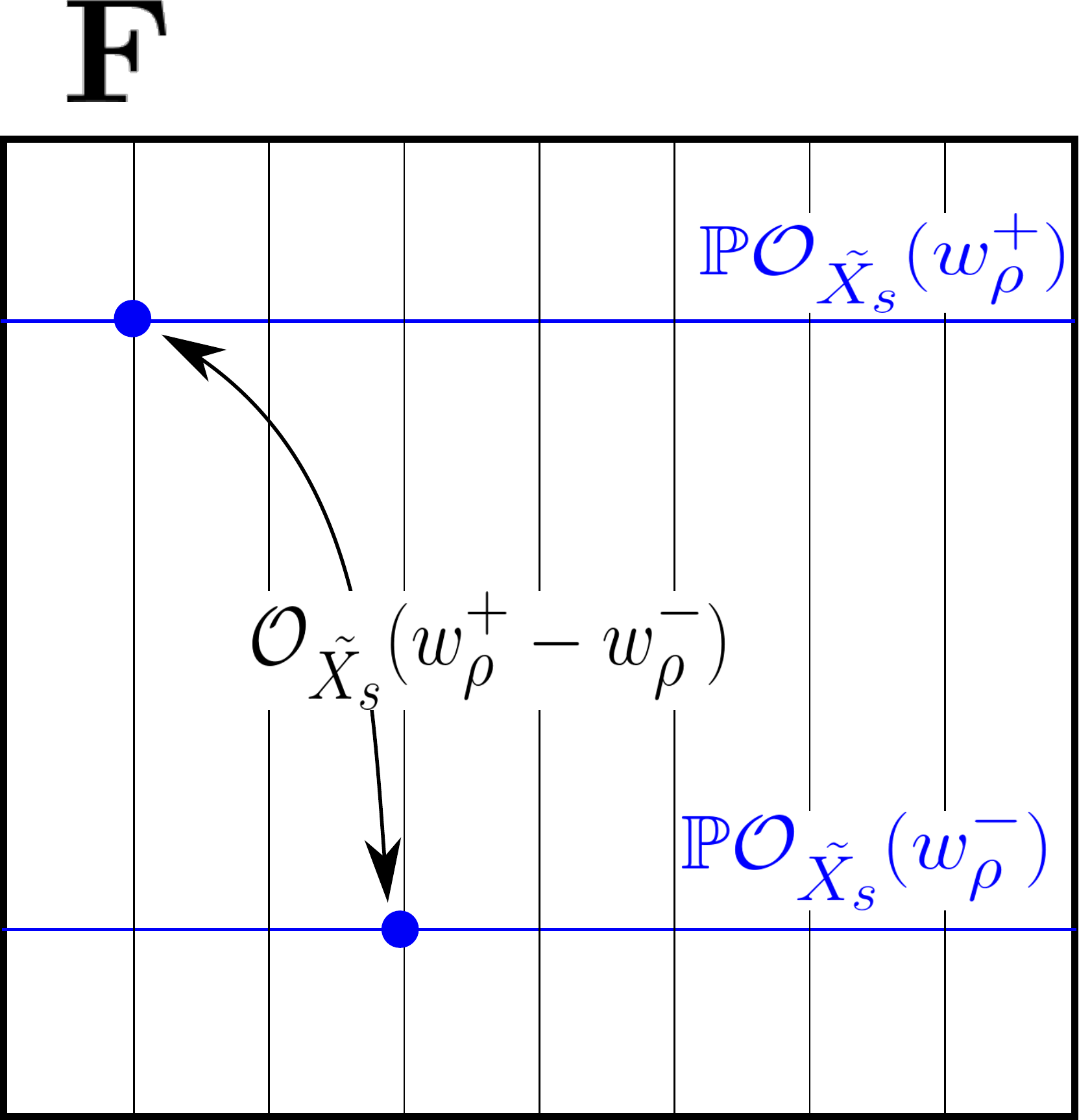}
\caption{Resolution of the compactified Jacobian.}
\label{p1bundle}
\end{figure}
\end{center}
\end{remark}

We will see that the singular Hitchin fiber ${\det}^{-1}(s)$, $s\neq 0$, is a union of two copies of ${\bf F}$. Before doing this, we introduce some notation. Let 
\[
\mathcal H = \mathcal H\setminus \mathcal N
\] 
be the complement of nilpotent cone and let $\mathcal H^{pairs}$ be the moduli space of pairs $(E, \theta)$ with $(E, {l}, \theta)$ in $\mathcal H$. Notice  that every Higgs bundle $(E, {l}, \theta)$ in $\mathcal H$ (and also in $\mathcal H^{pairs}$)  is irreducible, see \cite[Propositions 3.1 and 3.2]{FL22}.   We say that a pair $(E,\theta)$ is holomorphic at $\rho\in \{0,1,\lambda, t, \infty\}$ if ${\rm Res}(\theta,\rho)=0$. Each pair has at most one point with vanishing residual part, because the singular spectral curve has at most one singular (nodal) point.

\begin{lemma}\label{lemma:forgetful}
The forgetful map
\[
f: \mathcal H \to \mathcal H^{pairs}
\]
which forgets the parabolic structure, is the blowup at the locus ${\bf H}$ formed by pairs $(E, \theta)$ such that $\theta$ is holomorphic at some point $\rho\in \{0,1,\lambda, t, \infty\}$.  More precisely, $f$ is one to one in the complement $\mathcal H\setminus f^{-1}({\bf H})$ and $f^{-1}(E, \theta)$ is isomorphic to $\P^1$ for every $(E,\theta)\in {\bf H}$. 
\end{lemma}

\proof
If $\theta$ is nowhere-holomorphic, i.e., ${\rm Res}(\theta,\rho)\neq 0$ for every $\rho\in \{0,1,\lambda, t, \infty\}$, then the parabolic structure is determined by the  kernel of the residual part and the forgetful map is one to one. 

Now assume that ${\rm Res}(\theta,\rho)=0$ for some $\rho\in \{0,1,\lambda, t, \infty\}$ and we will show that  the fiber of the forgetful map is isomorphic to $\P^1$. Let 
\[
{l}({\rho}) = {l}\setminus \{l_{\rho}\}
\]
be the parabolic structure obtained by forgetting the direction over $\rho$ and let $(E, {l}({\rho}), \theta)$ be the corresponding Higgs bundle over $\P^1$ with four marked points 
\[
\Lambda_{\rho} = \{t_1, t_2, t_3, t_4\} = \{0,1,\lambda, t, \infty\}\setminus \{\rho\} .
\] 
The moduli space $Bun_{\mu}(0)$ parametrizing parabolic vector bundles $(E, {l}({\rho}))$ on $(\P^1, \Lambda_{\rho})$ of degree zero which are semistable with respect to weight 
$\mu = \left(\frac{1}{2}, \frac{1}{2}, \frac{1}{2}, \frac{1}{2} \right)$
is isomorphic to $\P^1$. A stable point of $Bun_{\mu}(0)$ has no automorphisms, besides trivial ones, then the fiber of $f$ is parametrized by the fifth parabolic direction $l_{\rho}\in \P E_{\rho}\simeq \P^1$, as we want. 

\begin{center}
\begin{figure}[h]
\centering
\includegraphics[height=1.3in]{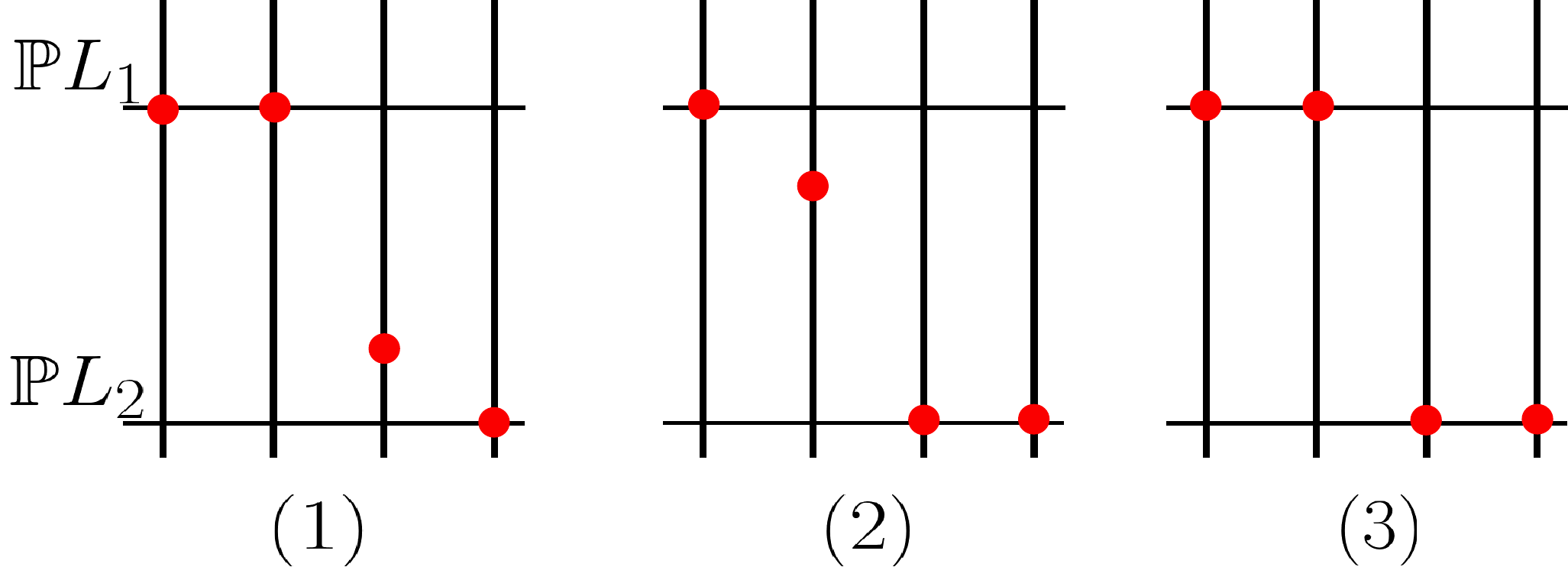}
\caption{Three $S$-equivalent parabolic structures giving a point in $Bun_{\mu}(0)$.}
\label{samestrict}
\end{figure}
\end{center}

It remains to consider strictly semistable points in  $Bun_{\mu}(0)$, there are exactly four of them, and each one is represented by three distinct quasi-parabolic structures giving the same $S$-equivalence class in $Bun_{\mu}(0)$, see Figure~\ref{samestrict}.  To see this remember that either $E=\mathcal O_{\P^1}\oplus \mathcal O_{\P^1}$ or  $E=\mathcal O_{\P^1}(1)\oplus \mathcal O_{\P^1}(-1)$. Since these four strictly semistable points are permuted by elementary transformations, we may assume that we are in one of  the three cases shown in Figure~\ref{samestrict}. Now, in the first two of them any Higgs field has vanishing determinant, and finally we arrive in the last case where $E=L_1\oplus L_2$, $L_i\simeq \mathcal O_{\P^1}$,
\[
l_1, l_2 \subset L_1 \;\;\text{and}\;\; l_3, l_4\subset L_2 
\]
and the Higgs field on $(E, {l}({\rho}))$ writes as 
\begin{eqnarray*}
\theta=\left(
\begin{array}{ccc} 
0 & a\frac{dx}{(x-t_1)(x-t_2)}  \\
b\frac{dx}{(x-t_3)(x-t_4)}   & 0 \\
\end{array}
\right).
\end{eqnarray*}
with  $a, b\in \C^*$. 
Adding the fifth parabolic direction $l_{\rho}$, if it does not lie in $L_1$ neither in $L_2$ then we may assume $l_{\rho} =  \begin{pmatrix}1\\ 1\end{pmatrix}$, $(E, {l})$ has no automorphisms and the fiber $f^{-1}(E,\theta)$ contains a $\C^*$ parametrized by
\begin{eqnarray*}
\theta_c=\left(
\begin{array}{ccc} 
0 & c a\frac{dx}{(x-t_1)(x-t_2)}  \\
c^{-1}b\frac{dx}{(x-t_3)(x-t_4)}   & 0 \\
\end{array}
\right)\;,\;\;c\in \C^*.
\end{eqnarray*}
It is worth noting that all the $\theta_c$ are equivalent in $\mathcal H^{pairs}$ because of the presence of automorphisms in $(E, {l}({\rho}))$, which are diagonal. To complete the fiber $f^{-1}(E,\theta)$ we have to add two points corresponding to either  $l_{\rho}\in L_1$ or $l_{\rho}\in L_2$. This finishes the proof of the lemma. 

\endproof

It follows from BNR correspondence \cite[Proposition 3.6]{BNR} that the fiber of the Hitchin map in the moduli  space $\mathcal H^{pairs}$ of pairs corresponds to the compactified  Jacobian variety $\overline{\rm Pic}^3(X_s)$ and the restriction of the forgetful map to ${\det}^{-1}(s)$ gives a map, still denoted by
\begin{eqnarray}\label{forgmap}
f:{\det}^{-1}(s) \to \overline{\rm Pic}^3(X_s).
\end{eqnarray}
To understand ${\det}^{-1}(s)$ we need the following result. 

\begin{lemma}\label{lemma:nodal}
Assume that the spectral curve $X_s$ has a  nodal singularity at $\rho\in \{0,1,\lambda, t, \infty\}$. There are bijective correspondences
\begin{itemize}
\item[(i)]\label{casenon} ${\rm Pic}^3(X_s) \leftrightarrow \left\{ (E, \theta)\in \mathcal H^{pairs}\;:\; \det\theta=s\;,\;\theta \;\text{is nowhere-holomorphic at}\;\rho\right\}$
\item[(ii)]\label{casehol} $\overline{\rm Pic}^3(X_s)\setminus  {\rm Pic}^3(X_s) \leftrightarrow \left\{ (E, \theta)\in \mathcal H^{pairs}\;:\; \det\theta=s\;,\;\theta \;\text{is holomorphic at}\;\rho\right\}$
\end{itemize}
\end{lemma}

\proof
The proof follows from  \cite[Proposition 3.5]{FL22}.

\endproof

In the case (i) of Lemma~\ref{lemma:nodal}, any Higgs field $\theta$ is {\it apparent } with respect to the parabolic direction over $\rho$, meaning that the parabolic direction $l_{\rho}$ is an eigendirection of the constant part of $\theta$. For instance, assuming that $l_{\rho} = \begin{pmatrix}1\\ 0\end{pmatrix}$ and $\rho=0$, we can write 
\begin{eqnarray}\label{appconst}
\theta=\left(
\begin{array}{ccc} 
ax & b \\
cx  & -ax \\
\end{array}
\right)\cdot \frac{dx}{x}
\end{eqnarray}
for suitable regular functions $a,b,c$ in a neighborhood of $\rho$, with $b(\rho)\neq 0$, because $\theta$ is nowhere-holomorphic at $\rho$. Since $X_s$ is singular over $\rho$ and then $\det \theta$ vanishes at order two, we conclude that $c(\rho)=0$, showing that $l_{\rho}$ is an eigendirection of the constant part of $\theta$. It is important to note that, after  an elementary transformation centered in $l_{\rho}$, the transformed Higgs field 
\begin{eqnarray}\label{appconst2}
\theta'=\left(
\begin{array}{ccc} 
ax & bx \\
c  & -ax \\
\end{array}
\right)\cdot \frac{dx}{x}
\end{eqnarray}
becomes holomorphic at $\rho$.  This discussion justifies the notation for ${\bf F}_{hol}$ and ${\bf F}_{app}$ in the next result.

\begin{center}
\begin{figure}[h]
\centering
\includegraphics[height=1.7in]{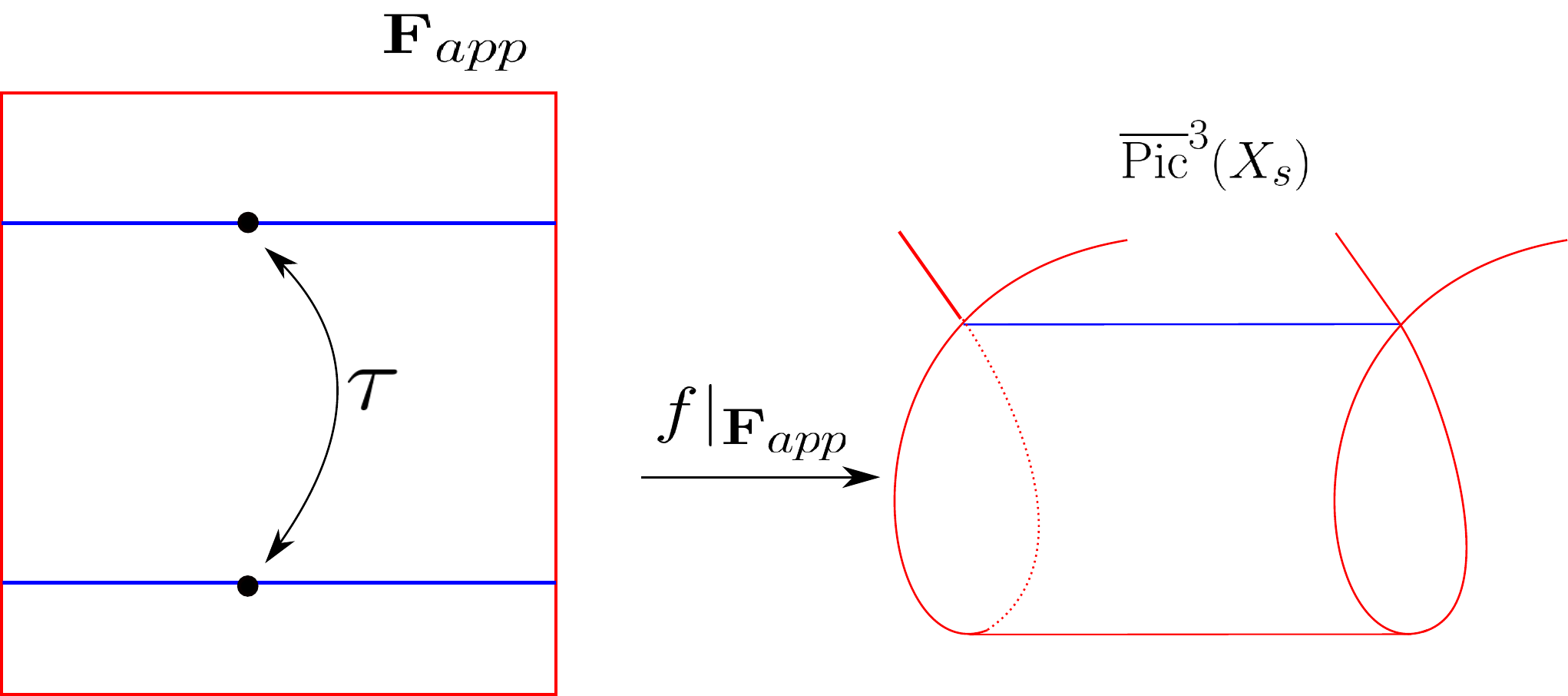}
\caption{Component ${\bf F}_{app}$.}
\label{Fapp}
\end{figure}
\end{center}

\begin{thm}\label{thm:sfiber}
Assume that the spectral curve $X_s$ has a  nodal singularity at $\rho\in \{0,1,\lambda, t, \infty\}$. The corresponding singular fiber ${\det}^{-1}(s)$ of the Hitchin map  has two irreducible components
\[
{\det}^{-1} (s) = {\bf F}_{hol}\cup{\bf F}_{app}
\]
which are isomorphic via any elementary transformation
\[
(elem_I)|_{{\bf F}_{hol}} : {\bf F}_{hol} \to {\bf F}_{app}
\]
where $I \subset \{0,1,\lambda, t, \infty\}$  contains $\rho$ and has even cardinality. Moreover:
\begin{enumerate}
\item Each component is a desingularization of $\overline{\rm Pic}^3(X_s)$, then isomorphic to ${\bf F}$, c.f. (\ref{compJac}), and the structure of $\P^1$-bundle in ${\bf F}_{hol}$ is given by 
\[
f|_{{\bf F}_{hol}}: {\bf F}_{hol}\to \tilde{X_s}\simeq \overline{\rm Pic}^3(X_s)\setminus  {\rm Pic}^3(X_s).
\] 
\item The map $f|_{{\bf F}_{app}}: {\bf F}_{app}\to  \overline{\rm Pic}^3(X_s)$ is a desingularization map. See Figure~\ref{Fapp}.
\item The intersection ${\bf F}_{hol}\cap{\bf F}_{app}$ is the union of the $0$-section and the $\infty$-section of ${\bf F}_{hol}$. See Figure~\ref{interComp}.
\end{enumerate}
\end{thm}

\proof
First, we identify $\overline{\rm Pic}^3(X_s)$ with the fiber of the Hitchin map 
\[
\det: \mathcal H^{pairs}\to \C^2
\] 
and ${\det}^{-1}(s)$ consists of $f^{-1}(\overline{\rm Pic}^3(X_s))$, where $f$ is the forgetful map (\ref{forgmap}). It follows from Lemmas~\ref{lemma:forgetful} and \ref{lemma:nodal} that ${\det}^{-1}(s)$ has two irreducible components, the strict transform of $\overline{\rm Pic}^3(X_s)$, which we call ${\bf F}_{app}$ and the $\P^1$-bundle ${\bf F}_{hol}$, which  is the blowup at the locus
\[
\left\{ (E, \theta)\in \mathcal H^{pairs}\;:\; \det\theta=s\;,\;\theta \;\text{is holomorphic at}\;\rho\right\}.
\]
This last is a copy of the elliptic curve $\tilde{X_s}$, because forgetting the parabolic direction over $\rho$ where $\theta$ is holomorphic, it can be identified with a fiber of the Hitchin map for moduli space of (irreducible) pairs $(E, \theta)$ over $\P^1$ with four parabolic points
\[
\{0,1,\lambda, t, \infty\}\setminus\{\rho\}.
\]
We conclude that ${\bf F}_{hol}$ is a $\P^1$-bundle over $\tilde{X_s}$. 

The elementary transformation $elem_I: \mathcal H \to \mathcal H$ is an isomorphism,  c.f (\ref{eleI}), and if $I$ contains $\rho$, $elem_I$ switches the  components ${\bf F}_{hol}$ and ${\bf F}_{app}$, see the discussion involving  (\ref{appconst}) and (\ref{appconst2}). In addition,  $f|_{{\bf F}_{app}}: {\bf F}_{app}\to \overline{\rm Pic}^3(X_s)$   is a birational morphism which is an isomorphism outside ${\bf F}_{hol} \cap {\bf F}_{app}$, and then 
$f|_{{\bf F}_{app}}$ is a desingularization map. 

We now study the intersection ${\bf F}_{hol} \cap {\bf F}_{app}$. Its restriction to  each fiber 
\[
f^{-1}(E,\theta)\simeq \P^1\subset {\bf F}_{hol} 
\]
corresponds to parabolic Higgs bundles $(E, {l}, \theta)$ with $\theta$ holomorphic at $\rho$ and apparent with respect to the parabolic direction $l_{\rho}$. Adding the fact that $X_s$ is nodal over $\rho$,  we can see that the constant part $\theta_{\rho}$ of $\theta$ has exactly two distinct eigendirections, it is an invertible matrix because otherwise $X_s$  would have a singularity of order bigger than two over $\rho$. In order to simplify notation, lets assume $\rho=t$, the other cases are similar. Any Higgs field in ${\bf F}_{hol}$ has determinant 
\[
\det \theta = s\cdot \frac{dx^{\otimes 2}}{x(x-1)(x-\lambda)}
\]
where $s\in \C^*$ is fixed and the constant part $\theta_t$ has determinant
\[
\det \theta_t =  \frac{s}{t(t-1)(t-\lambda)}
\]
which does not depend on $\theta$.  Therefore the intersection ${\bf F}_{hol} \cap {\bf F}_{app}$ is a union of two sections 
\[
\sigma_0, \sigma_{\infty}: \tilde{X_s}\to {\bf F}_{hol}
\]
where $\sigma_0$ is formed by eigendirections corresponding  to the eigenvalue $\sqrt{ \frac{-s}{t(t-1)(t-\lambda)}}$ and $\sigma_{\infty}$ corresponds to $-\sqrt{ \frac{-s}{t(t-1)(t-\lambda)}}$.

\endproof
%

\begin{remark}\rm
Via BNR correspondence, elements of ${\bf F}_{app}\backslash {\bf F}_{hol}$ correspond to line bundles on the nodal spectral curve $X_s$,  see Lemma~\ref{lemma:nodal} - (i). 
\end{remark}

We have seen that each irreducible component of ${\det}^{-1}(s)$  is a resolution of  $\overline{\rm Pic}^3(X_s)$. To recover $\overline{\rm Pic}^3(X_s)$ using ${\bf F}_{app}$ we must identify the $0$-section  and the $\infty$-section via the map 
\[
\tau:  \sigma_0(\tilde{X_s})\to  \sigma_{\infty}(\tilde{X_s})
\] 
which switches the two eigenvectors of the constant part of $\theta$. See Figure~\ref{Fapp}.

%

\begin{remark}\rm\label{rmk:multi}
Here,  we will see that $\tau$ consists of the translation $\mathcal O_{\tilde{X_s}}(w_{\rho}^+ - w_{\rho}^-)$, recovering Remark~\ref{rmk:Jac} from the modular point of view, in terms of elementary transformations on Higgs fields.

On the one hand, it is more convenient to work with ${\bf F}_{hol}$, instead of ${\bf F}_{app}$, because it contains a natural structure of $\P^1$-bundle given by the forgetful map $f$, and the resolution map is given by $f\circ elem_I: {\bf F}_{hol}\to \overline{\rm Pic}^3(X_s)$. On the other hand,  to recover the compactified Jacobian via ${\bf F}_{app}$, we need to identify the sections $\sigma_0$ and $\sigma_{\infty}$ gluing points in the same fiber of the forgetful map $f$, meaning that each point 
\[
(E,\theta)\in \overline{\rm Pic}^3(X_s)\setminus {\rm Pic}^3(X_s)\simeq {\tilde{X_s}}
\]
has exactly two representatives $\sigma_0(E,\theta)=(E, \theta, {l}_{\sigma_0})$ and  $\sigma_{\infty}(E,\theta)=(E, \theta, {l}_{\sigma_{\infty}})$ in ${\bf F}_{app}$, corresponding to the choices of eigendirections of the constant part of $\theta$. 

Coming back to ${\bf F}_{hol}$ using the involution $elem_I: {\bf F}_{hol}\to {\bf F}_{app}$, in order to obtain  $\overline{\rm Pic}^3(X_s)$ via  ${\bf F}_{hol}$ the $0$-section and the $\infty$-section must be identified via the map 
\begin{eqnarray}\label{secmap}
\iota := elem_I\circ \tau \circ elem_I: \sigma_0(\tilde{X_s}) \to \sigma_{\infty}(\tilde{X_s})
\end{eqnarray}
where $I$ has even cardinality and contains $\rho$. We will show  that this map corresponds to multiplication by $\mathcal O_{\tilde{X_s}}(w_{\rho}^+ - w_{\rho}^-)$. To do this,  let us first identify the elliptic curve $\tilde{X_s}$  with a fiber of the Hitchin map in the moduli space of pairs $(E, \theta)$ over $\P^1$ with four parabolic points $ \{0, 1, \lambda, t, \infty\}\setminus\{\rho\}$, and also with its Jacobian via BNR correspondence
\[
{\tilde{X_s}} \ni (E,\theta) \longleftrightarrow M_{\theta}\in {\rm Pic}(\tilde{X_s})\simeq \tilde{X_s}.
\]
There is a third identification for $\tilde{X_s}$,  for each $\theta$ with $\det \theta = s$, we identify $\tilde{X_s}$ with the curve of eigenvectors of $\theta$, and since $\theta$ is parabolic with respect to each one of the eigenvectors $w_{\rho}^+$ and $w_{\rho}^-$ of its  constant part  at $\rho$, then the variation of $M_{\theta}$ under an elementary transformation over $I$ centered in $w_{\rho}^{\pm}$, is given by \cite[Proposition 2.3]{FL22}.  Using this proposition, we see that the following diagram is commutative
 \begin{eqnarray*}
\xymatrix { 
\sigma_0(\tilde{X_s}) \ar@{<-}[d]_{\sigma_0} \ar@{->}[r]^{\iota}  & \sigma_{\infty}(\tilde{X_s}) \ar@{<-}[d]^{\sigma_{\infty}} \\
\tilde{X_s} \ar@{->}[r]_{\mathcal O_{\tilde{X_s}}(w_{\rho}^+ - w_{\rho}^-)}    &        \tilde{X_s}
}.
\end{eqnarray*}
\hfill$\lrcorner$
\end{remark}

\begin{remark}\rm\label{rkm:multi2} 
The structure of $\P^1$-bundle of ${\bf F}_{app}$ is obtained from ${\bf F}_{hol}$ via the isomorphism $elem_I: {\bf F}_{hol} \to {\bf F}_{app}$. Figure~\ref{interComp} shows a ruling of ${\bf F}_{app}$  intersecting  ${\bf F}_{hol}$. The whole Hitchin fiber ${\bf F}_{app}\cup {\bf F}_{hol}$ is a ``twisted product'' of an elliptic curve $\tilde{X_s}$ by a degenerate elliptic curve, meaning that a $\P^1$ of the ruling of ${\bf F}_{app}$ intersects two distinct $\P^1$'s of the ruling of ${\bf F}_{hol}$ and the intersection agrees with the multiplication by $\mathcal O_{\tilde{X_s}}(w_{\rho}^+ - w_{\rho}^-)$. The structure of this Hitchin fiber has been recently addressed by  C.T. Simpson in \cite[Discussion]{S18}, from the topological point of view.  
 \begin{center}
\begin{figure}[h]
\centering
\includegraphics[height=3.0in]{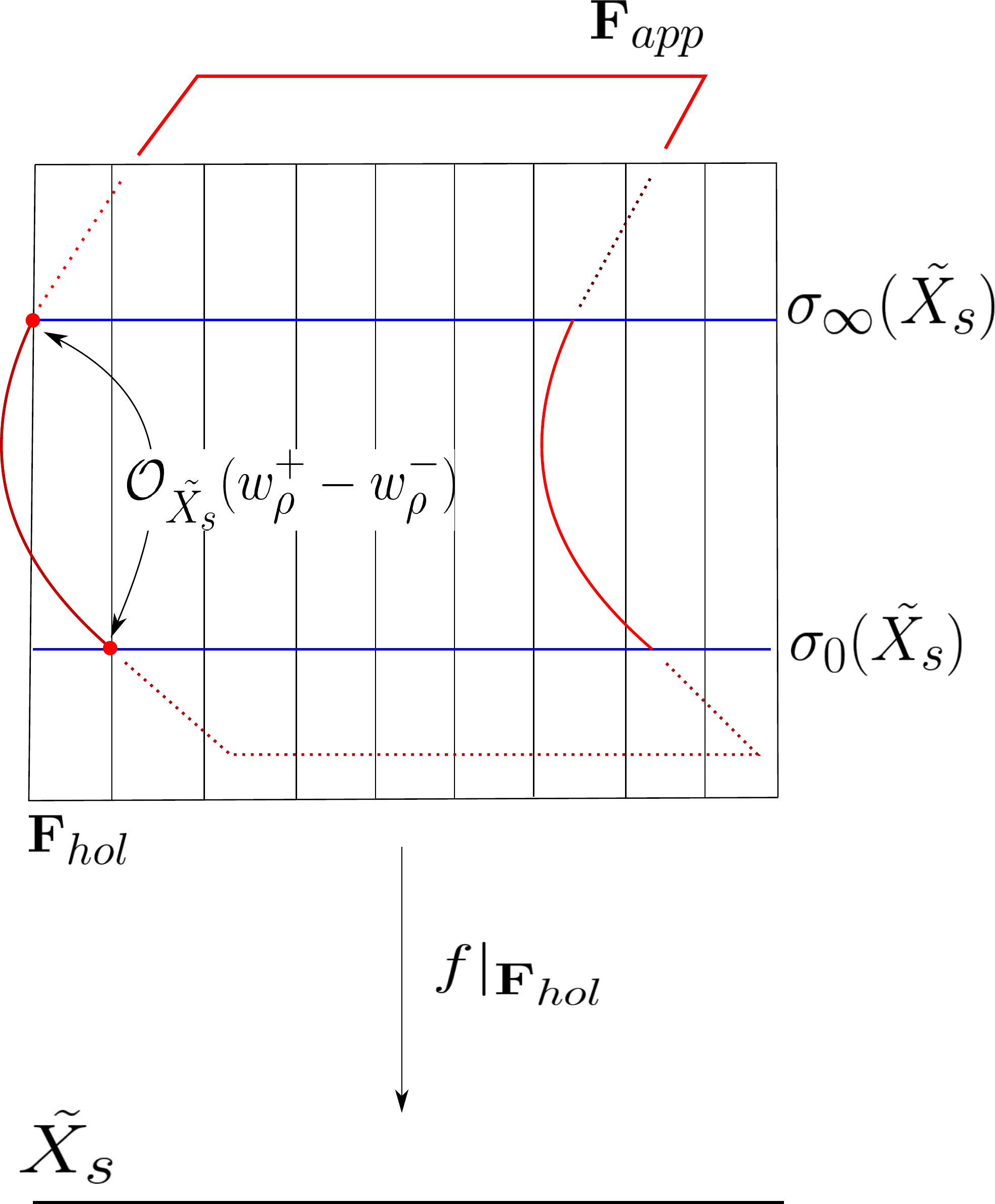}
\caption{Intersection of two components ${\bf F}_{hol}$ and ${\bf F}_{app}$.}
\label{interComp}
\end{figure}
\end{center}
\hfill$\lrcorner$
\end{remark}

\section{Connections}\label{section:connections}

Let $\mathcal C^{\nu}_n=\mathcal C^{\nu}(\P^1, \Lambda_n)$, $n\ge 5$, denote the moduli space of logarithmic connections over $\P^1$ of degree zero with polar divisor $\Lambda_n=t_1+\cdots+t_n$ supported on $n$ distinct points, and with prescribed eigenvalue vector $\nu = (\nu_1,\dots, \nu_n)\in\C^n$. An element of it consists of an isomorphism class $(E,\nabla)$, where $E$ is a rank two degree zero vector bundle over $\P^1$ endowed with a logarithmic connection, i.e. a $\mathbb C$-linear map
\[
\nabla\colon E \longrightarrow E\otimes \omega_{\P^1}(\Lambda_n)
\]
satisfying the Leibniz rule
\[
\nabla(a s)= s \otimes da + a \nabla (s)
\]
for (local) sections $s$ of $E$ and $a$ of $\mathcal{O}_{\P^1}$. In addition,  $\nabla$ is assumed to have vanishing trace and its 
 residue endomorphism ${\rm Res}_{t_i}(\nabla)$ over a given parabolic point $t_i$ has $\pm\nu_i$ as eigenvalues. 
 
We suppose that the eigenvalue vector $\nu$ is generic, meaning that $\nu_i\neq 0$, $\forall i$,  and  
\[
\sum\epsilon_i\nu_i\notin \mathbb Z
\]
for any choice of $\epsilon_i\in\{\pm 1\}$. 
From this, any connection is irreducible and the construction of the moduli space does not depend of a weight vector giving a stability notion.  The moduli space $\mathcal C^{\nu}_n$ is a smooth irreducible quasiprojective variety of dimension $2(n-3)$, see \cite{Inaba, IIS}. Note that 
each connection $\nabla$ on $E$ defines a unique parabolic structure,  by selecting the eigenspace $l_i\subset E\vert_{t_i}$ associated to $\nu_i$; 
therefore, $\mathcal C^{\nu}_n$ can equivalently be viewed as a moduli space of parabolic connections $(E, \nabla, {l})$. If a parabolic vector bundle $(E, {l})$ admits a connection like above, we say that it is {\it $\nu$-flat}. 

\subsection{Foliation conjecture} It follows from the work of Simpson \cite{S08} that there is a  decomposition of $\mathcal C^{\nu}_n$ obtained by looking at the limit of $c\cdot(E,\nabla, {l})$ as $c\to 0$. It turns out that for each weight vector $\mu$ and for $(E,\nabla,{l})\in\mathcal C^{\nu}_n$ there exists a unique limit
\[
(E, \theta, {\bf q}) = \lim_{c\to 0}  c\cdot(E, \nabla, {l}) \quad \in \mathcal H_{\mu}(\P^1, \Lambda_n, 0)
\] 
in the moduli space of $\mu$-semistable parabolic Higgs bundles, see also \cite[Proposition 4.1]{LSS}. This leads to an equivalence relation (depending on $\mu$) by assuming that two points of $\mathcal C^{\nu}_n$ are equivalent if their limits are the same.   We might equivalently consider the function 
\begin{eqnarray*}
\pi_{\mu}:\mathcal C^{\nu}_n &\to& \mathcal H_{\mu}(\P^1,\Lambda_n,0)\\
 (E, \nabla, {l}) &\mapsto& \lim_{c\to 0}  c \cdot(E, \nabla, {l})
\end{eqnarray*} 
and the decomposition of $\mathcal C^{\nu}_n$ given by fibers of $\pi_{\mu}$.  The {\it foliation conjecture} \cite[Question 7.4]{S08}, in this case, predicts that there is  a Lagrangian (regular) foliation $\mathcal F_{\mu}$ whose leaves are closed and  coincide with fibers of $\pi_{\mu}$. The Lagrangian property has already been proved by Simpson in \cite{S08}. The whole  conjecture has been proved  in \cite[Corollaries 5.7 and 6.2]{LSS} for the moduli space of connections over the four punctured projective line when the weight vector is generic,  and recently \cite{HHZ} deals with  the five punctured projective line by assuming  the weight vector $\mu$ satisfies $\sum\mu_i<1$, which lies in the unstable zone. For the unstable zone we mean the locus of weight vectors $\mu$ such that  any parabolic vector bundle is $\mu$-unstable. It is known that there is a polytope $\Delta\subset [0,1]^n$ consisting of weight vectors $\mu$ such that $Bun_{\mu}(\P^1,\Lambda_n,0)$ is nonempty \cite{B}, so the unstable zone consists of the complement of $\Delta$.

We will prove below that in the interior of $\Delta$ the foliation conjecture is sensitive to weight change, namely  it is true for the central weight $\mu_c=\left(\frac{1}{2},\frac{1}{2},\frac{1}{2},\frac{1}{2}, \frac{1}{2}\right)$ but it turns to be false if $\mu=\left(\frac{3}{4},\frac{1}{4},\frac{1}{4},\frac{1}{4}, \frac{1}{4}\right)$.  Even though the corresponding decompositions given by fibers of $\pi_{\mu_c}$ and $\pi_{\mu}$ share a Zariski open subset.   

\subsection{Foliation $\mathcal F_{Bun}$}
We shall consider the non-separated scheme $\mathcal P$ of rank two undecomposable parabolic vector  bundles over $(\P^1, \Lambda_n)$ and the corresponding forgetful map $\mathcal C^{\nu}_n\to \mathcal P$, sending $(E, \nabla, {l})$ to $(E, {l})$.  

\begin{prop}\label{prop:Fbun}
Each fiber of $\mathcal C^{\nu}_n\to \mathcal P$ is isomorphic to the affine space   $\C^n$ and they fit together into a regular foliation $\mathcal F_{Bun}$ on $\mathcal C^{\nu}_n$.  
\end{prop}

\proof
It follows from \cite[Proposition 3.1]{LS} that the following notions are equivalent 
\[
\nu\text{-flat} \Leftrightarrow \text{undecomposable} \Leftrightarrow \text{simple}
\] 
where simple means that any automorphism of $E$ preserving parabolic directions is scalar.  This implies that each fiber of  $\mathcal C^{\nu}_n\to \mathcal P$ is isomorphic to an affine space $\C^n$. Now, given $(E, \nabla, {l})$ in $\mathcal C^{\nu}_n$, by \cite[Proposition 3.4]{LS} the underlying parabolic vector bundle is $\mu$-stable for a convenient choice of weight vector $\mu$.  The local chart $Bun_{\mu}(\P^1,\Lambda,0)$ of $\mathcal P$ is a smooth irreducible projective variety and the restriction of  $\mathcal C^{\nu}_n\to \mathcal P$ to this chart gives a foliated neighborhood of  $(E, \nabla, {l})$ whose leaves coincide with fibers of  $\mathcal C^{\nu}_n\to \mathcal P$.  Varying the weight vector $\mu$ in all possible chambers, these foliated  neighborhoods fit together into a regular foliation $\mathcal F_{Bun}$ on $\mathcal C^{\nu}_n$. 

\endproof

The foliation $\mathcal F_{Bun}$ of Proposition~\ref{prop:Fbun} plays an important role when the weight vector is in the interior of the polytope $\Delta$.  In fact when $(E, {l})$ is $\mu$-stable, the limit  $\lim_{c\to 0}  c\cdot (E, \nabla, {l})$ consists of $(E, 0, {l})$,  then if $\mathcal U_{\mu}$ denotes the Zariski open subset formed by connections $(E, \nabla, {l})$ with $(E, {l})\in Bun_{\mu}(\P^1,\Lambda,0)$, the decomposition given by fibers of $\pi_{\mu}$ coincides with $\mathcal F_{Bun}$ when restricted to $\mathcal U_{\mu}$. In particular, we have the following result. 

\begin{prop}\label{prop:folcoic}
Assume that $\mu$ lies in the stable zone, i.e. it is in the interior of $\Delta$. If the foliation conjecture is true, that is, if the fibers of $\pi_{\mu}$ fit into a  regular foliation $\mathcal F_{\mu}$ then $\mathcal F_{\mu}=\mathcal F_{Bun}$. 
\end{prop}

\proof
By the discussion above, both foliations coincide on a nonempty open Zariski subset $\mathcal U_{\mu}$, then they must coincide everywhere. 

\endproof

\subsection{Variation with weights} In the next result we prove that given $n\ge 5$ there is a weight vector $\mu$ such that the  foliation conjecture \cite{S08} is false in the case $\P^1$ minus $n$ points. 

\begin{prop}\label{prop:foliation}
Let $n\ge 5$. The foliation conjecture in the moduli space $\mathcal C^{\nu}_n$ of logarithmic connections over the $n$ punctured projective line is false when $\mu=\left(\mu_1,\dots, \mu_n\right)$, $\mu_{n-2}=\frac{n-2}{n-1}$ and $\mu_i=\frac{1}{n-1}$, $\forall i\neq n-2$.
\end{prop}

\proof
Up to an automorphism of $\P^1$, we may assume $t_{n-2}=0$, $t_{n-1}=1$ and $t_n=\infty$.  By performing one elementary transformation over the parabolic point $t_{n-2}$, we go to the democratic weight $\mu' = \left(\frac{1}{n-1},\dots, \frac{1}{n-1}\right)$ and the determinant line bundle becomes odd. By \cite[Proposition 3.7]{LS},  the moduli space $Bun_{\mu'}(\P^1,\Lambda_n, -1)$ is isomorphic to $\P^{n-3}$, which gives the same conclusion to  $Bun_{\mu}(\P^1,\Lambda_n, 0)$. Fibers of $\pi_\mu$ over a point $(E, {l}, 0)$ with $(E, {l})$ in  $Bun_{\mu}(\P^1,\Lambda_n, 0)$ agree with leaves of $\mathcal F_{Bun}$. 

We now consider a $\nu$-flat parabolic vector bundle $(E, {l})$ which does not belong to $Bun_{\mu}(\P^1,\Lambda_n, 0)$ and let's investigate the fiber of $\pi_{\mu}$ over this point. Let us assume that $E=\mathcal O_{\P^1}\oplus \mathcal O_{\P^1}$, and parabolic directions are assumed to be
\[
l_{t_1}= \begin{pmatrix}u\\ 1\end{pmatrix}, l_{t_2}=\cdots=l_{t_{n-3}} =l_0= \begin{pmatrix}0\\ 1\end{pmatrix},  l_1= \begin{pmatrix}1\\ 1\end{pmatrix}, l_{\infty}= \begin{pmatrix}1\\ 0\end{pmatrix}.
\]
The parabolic structure is actually determined by $u\in \C$, so we denote by $(E, {l}_u)$ the corresponding parabolic vector bundle. Note that the embedding $\mathcal O_{\P^1}\to E$ corresponding to the second factor is a destabilizing subbundle, which makes $(E, {l}_u)$  $\mu$-unstable. Let us denote by $\C_u^{n-3}$ the space of connections over $(E, {l}_u)$. By \cite[Section 5.1]{LS}, this space is formed by connections $\nabla = \nabla_0+a_1\theta_1+\cdots +a_{n-3}\theta_{n-3}$, $(a_1,\dots,a_{n-3})\in \C_u^{n-3}$, where 
\begin{eqnarray*}
\nabla_0 = &d&+\left(
\begin{array}{ccc} 
-\nu_0 & 0 \\
\rho  & \nu_0 \\
\end{array}
\right) \frac{dx}{x}+
\left(
\begin{array}{ccc} 
-\nu_1-\rho & 2\nu_1+\rho \\
-\rho  & \nu_1+\rho \\
\end{array}
\right) \frac{dx}{x-1} 
+
\left(
\begin{array}{ccc} 
-\nu_{t_1} & 2\nu_{t_1} u \\
0  & \nu_{t_1} \\
\end{array}
\right) \frac{dx}{x-t_1} \\
&+&
\sum_{i=2}^{n-3}
\left(
\begin{array}{ccc} 
-\nu_{t_i} & 0 \\
0  & \nu_{t_i} \\
\end{array}
\right) \frac{dx}{x-t_i}\quad ,\text{with}\quad \rho=-\sum_{i=1}^{n-3} \nu_{t_i}-\nu_0-\nu_1-\nu_{\infty}  
\end{eqnarray*}
the Higgs fields are 
\begin{eqnarray*}
\Theta_1 = \left(
\begin{array}{ccc} 
0 & 0 \\
1-u  & 0 \\
\end{array}
\right) \frac{dx}{x}+
\left(
\begin{array}{ccc} 
u & -u \\
u  & -u \\
\end{array}
\right) \frac{dx}{x-1} +
\left(
\begin{array}{ccc} 
-u & u^2 \\
-1  & u \\
\end{array}
\right) \frac{dx}{x-t_1}  
\end{eqnarray*}
and
\begin{eqnarray*}
\Theta_i = \left(
\begin{array}{ccc} 
0 & 0 \\
1  & 0 \\
\end{array}
\right) \frac{dx}{x} +
\left(
\begin{array}{ccc} 
0 & 0 \\
-1  & 0 \\
\end{array}
\right) \frac{dx}{x-t_i} \quad , i=2,\dots,n-3 .
\end{eqnarray*}
Then we can take the gauge transformation rescaling by $c$ in the second component
\begin{eqnarray*}
g_c = \left(
\begin{array}{ccc} 
1 & 0 \\
0  & c \\
\end{array}
\right)
\end{eqnarray*}
to get 
\begin{eqnarray*}
\lim_{c\to 0} g_c(c\nabla)g_c^{-1} =\theta(a_1)= \left(
\begin{array}{ccc} 
0 & \beta \\
0  & 0 \\
\end{array}
\right)
\end{eqnarray*}
where 
\[
\beta = (2\nu_1+\rho-a_1u)\frac{dx}{x-1}+(2\nu_{t_1} + a_1u^2)\frac{dx}{x-t_1}.
\]
When $c$ goes to $0$, the parabolic structure projects to $(E, {\bf q})$ where 
\[
q_{t_1}= q_1 = q_{\infty}=\begin{pmatrix}1\\ 0\end{pmatrix}, l_{t_2}=\cdots=l_{t_{n-3}} =l_0= \begin{pmatrix}0\\ 1\end{pmatrix}
\]
and the limit Higgs bundle
\[
(E, \theta(a_1), {\bf q}) = \lim_{c\to 0}  c\cdot(E, \nabla, {l}_u) 
\] 
is stable with respect to the weight $\mu$, indeed the destabilizing subbundle  $\mathcal O_{\P^1}\to E$ given by the second factor is not invariant under $\theta(a_1)$. We are not able to eliminate  the parameter $a_1$ from $\theta(a_1)$ using automorphisms of $(E, {\bf q})$, so this computation shows that the leaf $\C^{n-3}_u$ of $\mathcal F_{Bun}$ is not contracted by $\pi_{\mu}$. This implies that fibers of $\pi_\mu$ and leaves of $\mathcal F_{Bun}$ do not agree everywhere. In view of Proposition~\ref{prop:folcoic},  we conclude that fibers of $\pi_\mu$ do not fit into a regular foliation on $\mathcal C^\nu_n$.  

\endproof

Our next result shows that, in the case $n=5$,  $\mathcal F_{Bun}$ can be realized as fibers of $\pi_{\mu_c}$, for  the central weight $\mu_c = \left(\frac{1}{2},\frac{1}{2},\frac{1}{2},\frac{1}{2},\frac{1}{2}\right)$.

\begin{thm}\label{thm:foliation}
 The foliation conjecture in the moduli space $\mathcal C^{\nu}_5$ of logarithmic connections over the five punctured sphere is true when  $\mu_c=\left(\frac{1}{2},\frac{1}{2},\frac{1}{2},\frac{1}{2}, \frac{1}{2}\right)$.
\end{thm}

\proof
Let $(E, \nabla, {l})$ be an element of $\mathcal C^{\nu}_5$. It follows from \cite[Corollary 3.3]{LS} that either $E = \mathcal O_{\P^1}(1)\oplus \mathcal O_{\P^1}(-1)$ or $E = \mathcal O_{\P^1}\oplus \mathcal O_{\P^1}$.

The locus of fixed points by the $\C^*$-action on the moduli space $\mathcal H$ of Higgs bundles is the union of $\mathcal S$, corresponding to $(E, 0 , {l})$ with $(E,{l})$ $\mu$-stable, and the $16$ Hodge bundles $\Theta_i$, see Definition~\ref{def:Hodge} and Theorem~\ref{thm:nilpP1}. A fiber of $\pi_{\mu_c}$ over a point $(E, 0 , {l})$ consists of a leaf of $\mathcal F_{Bun}$,  so it remains to consider the other $16$ points. 

We will show that there are exactly 16 $\mu$-unstable $\nu$-flat parabolic vector bundles. Assuming that $(E,{l})$ is $\mu_c$-unstable, there exists a destabilizing subbundle $L\subset E$ satisfying 
\[
-2\deg L - \frac{m}{2}+\frac{5-m}{2}<0
\]
where $m$ is the number of parabolic directions lying in $L$. This gives $\deg L \in \{-1,0,1\}$. If $\deg L = 1$ the there is at least one parabolic direction in $L$ and at least two parabolic directions outside $L$, otherwise $(E, {l})$ would be decomposable. Up to performing an elementary transformation over these two parabolic directions outside $L$, we may assume that $L$  has degree zero and  $E = \mathcal O_{\P^1}\oplus \mathcal O_{\P^1}$. The same reasoning can be applied to the case where $L$ has degree $-1$, indeed here all parabolic directions must lie in $L$,  we then apply an elementary transformation over two of them. Therefore we may assume that $L$ has degree zero, $E = \mathcal O_{\P^1}\oplus \mathcal O_{\P^1}$ and exactly three parabolic directions lie in $L$ (more than three implies $(E, {l})$ undecomposable). We then arrive, up to elementary transformations, in the following case:  $E = \mathcal O_{\P^1}\oplus \mathcal O_{\P^1}$ and 
\begin{eqnarray}\label{parabb}
l_{0}= l_{\lambda} = l_{t}= \begin{pmatrix}0\\ 1\end{pmatrix}, l_1= \begin{pmatrix}1\\ 1\end{pmatrix},   l_{\infty}= \begin{pmatrix}1\\ 0\end{pmatrix}.
\end{eqnarray}
This implies that there are exactly $16$ $\mu$-unstable $\nu$-flat parabolic vector bundles, they are in the same orbit of the group ${\bf El}$ of elementary transformations.  

The space of connections over the parabolic bundle (\ref{parabb}) is formed by $\nabla = \nabla_0+a_1\theta_1+a_2\theta_2$, $a_1,a_2\in\C$, where   
\begin{eqnarray*}
\nabla_0 = &d&+\left(
\begin{array}{ccc} 
-\nu_0 & 0 \\
\rho  & \nu_0 \\
\end{array}
\right) \frac{dx}{x}+
\left(
\begin{array}{ccc} 
-\nu_1-\rho & 2\nu_1+\rho \\
-\rho  & \nu_1+\rho \\
\end{array}
\right) \frac{dx}{x-1} \\
&+&
\left(
\begin{array}{ccc} 
-\nu_{\lambda} & 0 \\
0  & \nu_{\lambda} \\
\end{array}
\right) \frac{dx}{x-\lambda} +
\left(
\begin{array}{ccc} 
-\nu_{t} & 0 \\
0  & \nu_{t} \\
\end{array}
\right) \frac{dx}{x-t}\quad ,\text{with}\quad \rho=-\sum \nu_{i}  
\end{eqnarray*}
the Higgs fields are 
\begin{eqnarray*}
\theta_1 = \left(
\begin{array}{ccc} 
0 & 0 \\
1  & 0 \\
\end{array}
\right) \frac{dx}{x}+
\left(
\begin{array}{ccc} 
0 & 0 \\
-1  & 0 \\
\end{array}
\right) \frac{dx}{x-\lambda}  
\end{eqnarray*}
and
\begin{eqnarray*}
\theta_2 = \left(
\begin{array}{ccc} 
0 & 0 \\
1  & 0 \\
\end{array}
\right) \frac{dx}{x} +
\left(
\begin{array}{ccc} 
0 & 0 \\
-1  & 0 \\
\end{array}
\right) \frac{dx}{x-t}.
\end{eqnarray*}

Then  the gauge transformation
\begin{eqnarray*}
g_c = \left(
\begin{array}{ccc} 
1 & 0 \\
0  & c \\
\end{array}
\right)
\end{eqnarray*}
gives
\begin{eqnarray*}
\lim_{c\to 0} g_c(c\nabla)g_c^{-1} =\theta= \left(
\begin{array}{ccc} 
0 & \beta \\
0  & 0 \\
\end{array}
\right)
\end{eqnarray*}
where 
\[
\beta = (2\nu_1+\rho)\frac{dx}{x-1}.
\]
Note that using an automorphism of the projected parabolic vector bundle, we can eliminate the constant $2\nu_1+\rho$ from $\beta$.  Indeed when $c$ goes to $0$, the parabolic structure projects to $(E, {\bf q})$ where 
\[
q_{0}= q_{\lambda} = q_{t}=\begin{pmatrix}0\\ 1\end{pmatrix}, l_{1}= l_{\infty}=\begin{pmatrix}1\\ 0\end{pmatrix}
\]
The conclusion is that the limit Higgs bundles $ \lim_{c\to 0}  c\cdot(E, \nabla, {l}) $ is one of the $16$ Hodge bundles. 
Therefore, any fiber of $\pi_{\mu_c}$ coincides with a leaf of $\mathcal F_{Bun}$. 
\endproof


\begin{thebibliography}{9}

%
%
%

\bibitem{AFKM21}
C. Araujo, T. Fassarella, I. Kaur, and A. Massarenti, 
{\it On automorphisms of moduli spaces of parabolic vector bundles}, 
International Mathematics Research Notices, {\bf 3} (2021) 2261-2283.

\bibitem{B}
S. Bauer, 
{\it Parabolic bundles, elliptic surfaces and $SU(2)$-representation spaces of genus zero Fuchsian groups.}
Math. Ann. {\bf 290} (1991) 509-526.


\bibitem{BNR}
A. Beauville, M.S. Narasimhan and S. Ramanan, 
{\it Spectral curves and the generalised theta divisor.}
J. reine angew. Math. (Crelles Journal), {\bf 398} (1989) 169-179.



%
%
%
\bibitem{CHM12}
 M. A. de Cataldo, T. Hausel, L. Migliorini. 
 {\it Topology of Hitchin systems and Hodge theory of character varieties: the case A1.} 
 Ann. of Math. {\bf 175} (2012) 1329-1407.

\bibitem{FL}
T. Fassarella and F. Loray,
{\it Flat parabolic vector bundles on elliptic curves.}
J. reine angew. Math. (Crelles Journal), {\bf 761} (2020) 81-122.

\bibitem{FL22}
T. Fassarella and F. Loray,
{\it Hitchin fibration under ramified coverings}
(2022)  	arXiv:2208.10130. 


%
%
%
%

\bibitem{GO13}
P. B. Gothen and A. G. Oliveira, 
{\it The singular fiber of the Hitchin map}, 
Int. Math. Res. Not. IMRN, {\bf 5} (2013) 1079-1121. 

\bibitem{HT}
T. Hausel and M. Thaddeus,
{\it Mirror symmetry, Langlands duality, and the Hitchin system},
Invent. Math., {\bf 153} (2003) 197-229. 

\bibitem{H87}
N. J. Hitchin, 
{\it The self-duality equations on a Riemann surface}, 
Proc. London Math. Soc. {\bf 3} 55(1) (1987) 59-126.



\bibitem{Hor}
J. Horn, 
{\it Semi-abelian spectral data for singular fibres of the SL(2,C)-Hitchin system}, 
Int. Math. Res. Not. IMRN {\bf 10} (2020) 3860-3917.

\bibitem{Hor2}
J. Horn,
{\it $\frak{sl}_2$-Type singular fibres of the symplectic and odd orthogonal Hitchin system},
Journal of Topology {\bf 15} (2022) 1-38. 


\bibitem{HHZ}
Z. Hu, P. Huang and R. Zong,
{\it Moduli spaces of parabolic bundles over $\P^1$ with five marked points},
(2021), arXiv:2108.08994. 

\bibitem{Inaba}
M.-a. Inaba,
{\it Moduli of parabolic connections on curves and the Riemann-Hilbert	correspondence}, 
J. Algebraic Geom. 22(3) (2013) 407-480.


\bibitem{IIS}
M.-a. Inaba, K.~Iwasaki, and M.-H. Saito.
{\it Moduli of stable parabolic connections, Riemann-Hilbert
	correspondence and geometry of Painlev\'{e} equation of type {VI}. {I}.}
Publ. Res. Inst. Math. Sci. 42(4) (2006) 987-1089.


\bibitem{LS}
F. Loray and M.-H. Saito,
{\it Lagrangian fibrations in duality on moduli spaces of rank 2 logarithmic connections over the projective line.} 
Int. Math. Res. Not. IMRN, {\bf 4} (2015) 995-1043. 

\bibitem{LSS}
F. Loray, M.-H. Saito and C. Simpson,
{\it Foliations on the moduli space of rank two connections on the projective line minus four points},
S\'eminaires et Congr\`es, {\bf 27} (2013) 115-168.


%
%
%
%
%


\bibitem{Ngo1}
B. C. Ng\^o, 
{\it Fibration de Hitchin et endoscopie}, 
Invent. Math. {\bf 164} (2006) 399-453.

\bibitem{OS79}
T. Oda and C. S. Seshadri,
{\it Compactifications of the generalized Jacobian variety,}
Trans. Amer. Math. Soc. {\bf 253} (1979) 1-90. 

\bibitem{S08}
C. T. Simpson. 
{\it Iterated destabilizing modifications for vector bundles with connection.} 
Vector Bundles and Complex Geometry (Conference in Honor of Ramanan, Miraflores 2008)
Contemporary Math. 522 (2010) 183-206.

\bibitem{S18}
C. T. Simpson, 
{\it An explicit view of the Hitchin fibration on the Betti side for $\P^1$ minus five points.} 
Geometry and Physics: Volume 2: A Festschrift in honour of Nigel Hitchin (2018): 705.

%
%
%
\bibitem{Sim90}
C. T. Simpson, 
{\it Harmonic bundles on noncompact curves}, 
J. Amer. Math. Soc. 3(3) (1990) 713-770.

\bibitem{Sim92}
C. T. Simpson, 
{\it Higgs bundles and local systems}, 
Inst. Hautes  \'Etudes Sci. Publ. Math. {\bf 75} (1992) 5-95.

%
%


\bibitem{Yo93}
K. Yokogawa
{\it Compactification of moduli of parabolic sheaves and moduli of parabolic Higgs sheaves},
J. Math. Kyoto Univ. {\bf 33} (1993) 451-504. 

\bibitem{Yo95}
K. Yokogawa,
{\it Infinitesimal deformation of parabolic Higgs sheaves},
International Journal of Mathematics, {\bf 6} (1995) 125-148.
 
\end{thebibliography}
\end{document}